\newif\ifpreprint
\preprinttrue

\ifpreprint
\documentclass[12pt]{article}

\usepackage{fullpage}

\let\oldparagraph\paragraph

\renewcommand{\paragraph}[1]{\oldparagraph{#1.}}

\usepackage{amsmath,amssymb,amsthm,amsfonts}
\usepackage{xcolor}

\usepackage[noend]{algpseudocode} 

\usepackage{csvsimple}	\usepackage{booktabs} \let\oldtabular\tabular
\let\endoldtabular\endtabular
\renewenvironment{tabular}{\small\oldtabular}{\endoldtabular}
\usepackage{adjustbox}  

\usepackage[colorlinks=true, citecolor=cyan, linkcolor=cyan, urlcolor=cyan]{hyperref}

\usepackage{caption}
\newtheorem{theorem}{Theorem}
\newtheorem{lemma}{Lemma}
\newtheorem{corollary}{Corollary}[theorem]  \newtheorem{proposition}{Proposition}
\newtheorem{assumption}{Assumption}

\theoremstyle{definition}
\newtheorem{definition}{Definition}

\theoremstyle{remark}
\newtheorem{remark}{Remark}
\newtheorem{example}{Example}
\newtheorem*{example*}{Example}

\allowdisplaybreaks

\usepackage[numbers,sort&compress]{natbib}
\else
\documentclass[twoside,11pt]{article}

\usepackage{amsmath,amsthm,amsfonts}
\usepackage[preprint]{jmlr2e}

\renewenvironment{proof}[1][Proof]{\par\noindent{\bf #1.\ }}{\hfill\BlackBox\par\addvspace{2mm}}

\newtheorem{assumption}[theorem]{Assumption}
\newtheorem*{example*}{Example}

\usepackage{booktabs}

\allowdisplaybreaks

\let\oldparagraph\paragraph
\renewcommand{\paragraph}[1]{\oldparagraph{#1.}}

\fi 

\usepackage{graphicx}
\usepackage{tikz-cd}
\usepackage{algorithm,algorithmicx}
\usepackage{multirow}
\usepackage{csvsimple}
\usepackage{enumitem,mathtools}
\usepackage{hyperref}
\hypersetup{hypertexnames=false}

\newcommand{\eg}{\textit{e.g.}}
\newcommand{\ie}{\textit{i.e.}}

\newcommand{\norm}[2][]{\left\|#2\right\|_{#1}}
\newcommand{\inner}[3][]{\left\langle#2,#3\right\rangle_{#1}}

\newcommand{\ones}{{\mathbf 1}}
\newcommand{\reals}{{\mathbf R}}
\newcommand{\symm}{{\mathbf S}}
\newcommand{\soc}{{\mathcal{K}_{\text{SOC}}}}
\newcommand{\prob}{{\mathbf P}}
\newcommand{\identity}{I}
\newcommand{\Tr}{\mathop{\bf tr}}

\newcommand{\Span}{\mathop{\bf span}}
\newcommand{\lambdamax}{{\lambda_{\rm max}}}

\newcommand{\Expect}{\mathop{\bf E{}}}
\newcommand{\Prob}{\mathop{\bf Prob}}
\newcommand{\argmin}{\mathop{\rm argmin}}

\DeclareMathOperator*{\esssup}{{ess}\,{sup}}
\DeclareMathOperator{\prox}{prox}

\newcommand{\lagrangian}{\mathbf{L}}

\newcommand{\var}{\mathbf{VaR}}
\newcommand{\cvar}{\mathbf{CVaR}}

\newcommand{\supp}{\mathrm{supp}}

\newcommand{\ambiset}[1][]{\mathcal{B}_{#1}}

\newcommand{\hG}{\widehat{G}}
\newcommand{\hf}{\hat{f}}
\newcommand{\hF}{\widehat{F}}
\newcommand{\hx}{\hat{x}}

\newcommand{\hZ}{\widehat{Z}}
\newcommand{\perfmetr}[1][]{\phi^{#1}}
\newcommand{\xclass}{\mathcal{X}}
\newcommand{\fclass}{\mathcal{F}}
\newcommand{\fclassX}{\mathcal{F}_\mathcal{X}}
\newcommand{\fclassprox}{\mathcal{G}}
\newcommand{\lclass}{\mathcal{L}}
\newcommand{\linear}{\mathcal{S}}
\newcommand{\domain}{\mathcal{D}}

\newcommand{\Aobj}{A_{\mathrm{obj}}}
\newcommand{\bobj}{b_{\mathrm{obj}}}
\newcommand{\cobj}{c_{\mathrm{obj}}}

\newcommand{\algo}{\mathcal{A}}
\newcommand{\probQ}{\mathbf{Q}}
\newcommand{\hprob}{\widehat{\prob}}
\newcommand{\measure}{\mathcal{M}}
  
\newcommand{\suppref}[1]{Appendix~\ref{#1}}

\newcommand{\myabstract}{We consider the problem of analyzing the probabilistic performance of first-order methods when solving convex optimization problems drawn from an unknown distribution only accessible through samples. By combining performance estimation and Wasserstein distributionally robust optimization, we formulate the analysis as a tractable conic program. Our approach unifies worst-case and average-case analyses by incorporating data-driven information from the observed convergence of first-order methods on a limited number of problem instances. This yields probabilistic, data-driven performance guarantees in terms of the expectation or conditional value-at-risk of the selected performance metric.
Our open-source implementation computes these guarantees directly from sampled algorithm trajectories using off-the-shelf conic solvers.
Experiments on convex quadratic minimization, real-data logistic regression using a credit-scoring dataset, and Lasso show that our method significantly reduces the conservatism of classical worst-case bounds and narrows the gap between theoretical and empirical performance.
}

\newcommand{\myack}{Bartolomeo Stellato and Vinit Ranjan are supported by the NSF CAREER Award ECCS-2239771 and the ONR YIP Award N000142512147. Jisun Park is supported by the National Research Foundation of Korea (NRF) grant funded by the Korea government (MSIT) (RS-2024-00353014) and the ONR YIP Award N000142512147. The authors are pleased to acknowledge that the work reported on in this paper was substantially performed using Princeton University’s Research Computing resources.}

\ifpreprint 

\title
{Data-Driven Analysis of First-Order Methods\\
via Distributionally Robust Optimization}

\author{
Jisun Park, Vinit Ranjan, and Bartolomeo Stellato
}
\else
\title{Data-Driven Analysis of First-Order Methods\\ via Distributionally Robust Optimization}

\author{\name Jisun Park \email jisunpark@princeton.edu \\
        \addr Department of Operations Research and Financial Engineering, Princeton University, Princeton, NJ 08544, USA\\
        Research Institute of Mathematics, Seoul National University, Seoul 08826, South Korea
        \AND
        \name Vinit Ranjan \email vranjan1@mit.edu \\
        \addr MIT Sloan School of Management, Cambridge, MA 02142, USA
        \AND
        \name Bartolomeo Stellato \email bstellato@princeton.edu \\
        \addr Department of Operations Research and Financial Engineering, Princeton University, Princeton, NJ 08544, USA}

\ShortHeadings{Data-Driven Analysis of First-Order Methods}{Park, Ranjan, and Stellato}
\firstpageno{1}
\fi 

\begin{document}

\maketitle

\ifpreprint
\begin{abstract}
    \myabstract
\end{abstract}
\else
\begin{abstract}
    \myabstract
\end{abstract}

\fi

\section{Introduction}
\label{sec:intro}
First-order methods have gained significant attention over the past decade due to their low computational cost per iteration, modest memory requirements, and ability to warm-start from previous solutions~\citep{beckFirstOrderMethodsOptimization2017,ryuLargeScaleConvexOptimization2022}.
These advantages make them valuable for large-scale problems in data science and machine learning, as well as embedded optimization applications in engineering and optimal control.

These algorithms date to the 1950s~\citep{MargueriteWolfe1956} and have been extended to nondifferentiable functions and constraints via proximal operators~\citep{parikhProximalAlgorithms2014} and to composite objectives via operator-splitting methods, including the alternating direction method of multipliers (ADMM)~\citep{boydstephenDistributedOptimizationStatistical2011} and the primal-dual hybrid gradient method (PDHG)~\citep{chambolleFirstOrderPrimalDualAlgorithm2011}.
Recently, several general-purpose solvers based on first-order methods have appeared, including PDLP~\citep{applegatePracticalLargeScaleLinear2021,applegateFasterFirstorderPrimaldual2023} for linear programs, OSQP~\citep{stellatoOSQPOperatorSplitting2020} for quadratic programs, and SCS~\citep{odonoghue:21,ocpb:16} and COSMO~\citep{garstkaCOSMOConicOperator2021} for semidefinite programs (SDPs).
Despite these advances, the behavior of these algorithms is still highly sensitive to data, often exhibiting slow convergence for badly scaled problems.

To design faster and more reliable first-order methods, a central question is to understand the accuracy of the solutions returned.
In convergence analyses, this translates into estimating the worst-case accuracy after a certain number of iterations over a given function class.
While most of these results focus on asymptotic convergence rates~\citep{bertsekasConvexOptimizationTheory2009,nesterovLecturesConvexOptimization2018}, a recent line of work called the performance estimation problem (PEP)~\citep{droriPerformanceFirstorderMethods2014,taylorExactWorstCasePerformance2017} focuses on constructing finite-step guarantees by solving a semidefinite program.
Interpolation conditions represent many convex function classes, while linear equations describing the algorithm represent virtually all fixed-step first-order methods used in practice.
This approach has led to several new results in the optimization community, including tighter convergence-rate proofs~\citep{taylorStochasticFirstorderMethods2019,upadhyayaAutomatedTightLyapunov2025} and faster optimization algorithms~\citep{jangComputerAssistedDesignAccelerated2025,kimOptimizedFirstorderMethods2016,parkExactOptimalAccelerated2022}.
However, the practical convergence behavior of these schemes can be very different from the worst-case bounds.
In some cases, first-order methods might converge significantly faster, as the following example illustrates.
\begin{example} Consider a least-squares problem with box constraints:
    \begin{equation}\label{prob:box_lstsq}
        \begin{array}{ll}
            \text{minimize} & f(x) = (1/2) \norm[2]{Ax - b}^2 \\
            \text{subject to} & 0 \le x \le 1,
        \end{array}
    \end{equation}
    where~$x \in \reals^d$ is the optimization variable, and~$A \in \reals^{n \times d}$ and~$b \in \reals^n$, with~$n<d$, are problem data.
    To solve~\eqref{prob:box_lstsq}, we apply projected gradient descent with step size~$1/L$, where~$L$ is the maximum eigenvalue of~$A^T A$:
    \begin{equation*}
        x^{k+1}
        = \Pi_{[0, 1]} \big( x^k - (1/L) \nabla f(x^k) \big)
        = \Pi_{[0, 1]} \big( x^k - (1/L) A^T(Ax^k - b) \big),
        \quad k=0,1,\dots,
    \end{equation*}
    where~$\Pi_{[0, 1]}$ is an elementwise projection onto the interval~$[0, 1]$.
    It is well known~\citep{beckFastIterativeShrinkageThresholding2009,parikhProximalAlgorithms2014} that the worst-case function-value gap is
    \begin{equation*}
        f(x^K) - f(x^\star) \le \frac{L}{2K} \| x^0 - x^\star\|^2,
    \end{equation*}
    where~$x^\star$ is an optimal solution of~\eqref{prob:box_lstsq}.
    However, if we specifically consider the problem instances with varying $b$, the actual performance can be significantly better.
    To see this, we construct the data matrix~$A = U \Lambda V^T$, with random orthogonal matrices~$U \in \reals^{n \times n}$ and~$V \in \reals^{d \times d}$,
    along with a rectangular matrix~$\Lambda \in \reals^{n \times d}$ whose only nonzero entries are the diagonal entries~$\{\Lambda_{ii} \}_{i=1}^n$ drawn from the uniform distribution $\mathcal{U}[0, \sqrt{L}]$.
    We create $1000$ instances of~\eqref{prob:box_lstsq} by sampling the corresponding vectors $b$ with i.i.d.\ entries $b_i \sim \mathcal{N}(0, 1)$.
    \begin{figure}[t]
        \centering
\includegraphics[width=0.8\textwidth]{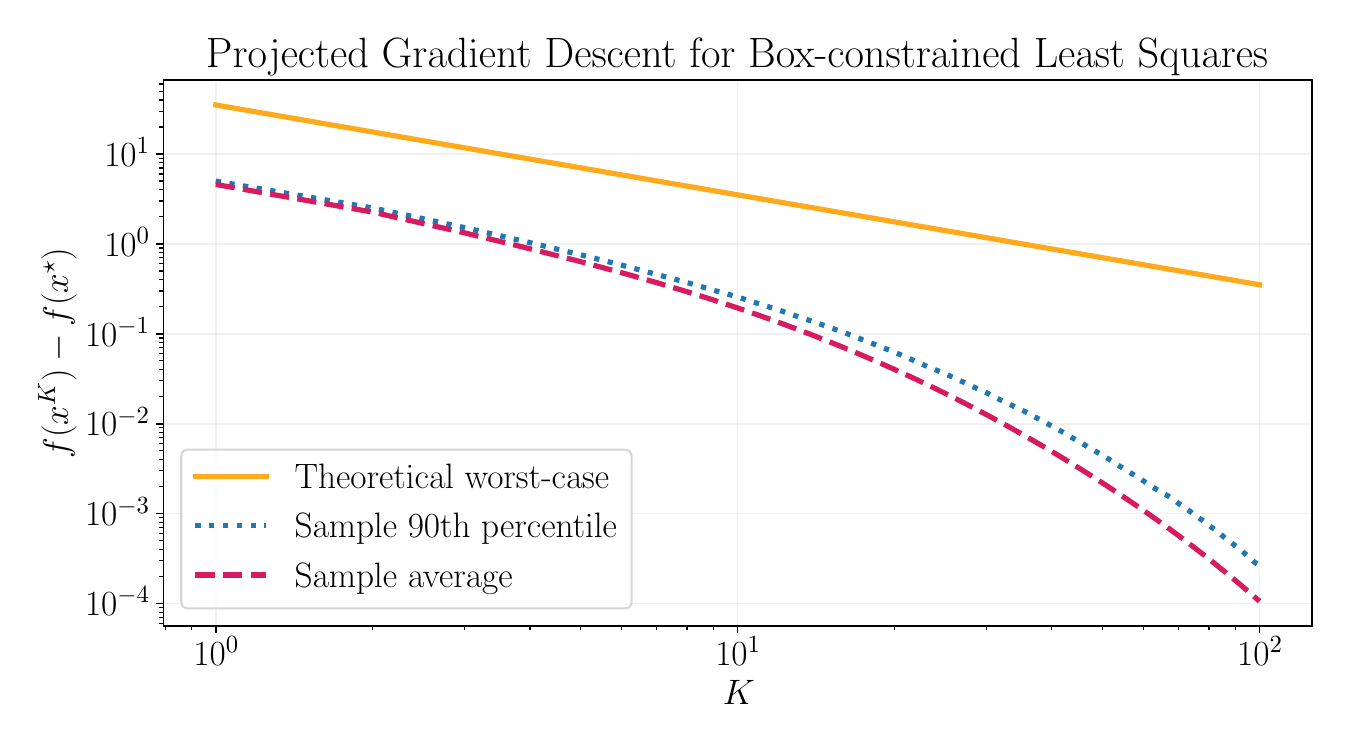}
        \caption{
            Function-value gap of projected gradient descent on the box-constrained least-squares problem~\eqref{prob:box_lstsq} with~$L = 1$.
            The curves show the theoretical worst-case bound, empirical mean, and empirical 90\textsuperscript{th} percentile across instances.
            The theoretical $O(K^{-1})$ worst-case bound on~$f(x^K) - f(x^\star)$ is pessimistic relative to the empirical curves.
        }
        \label{fig:box_lstsq_experiment}
    \end{figure}
    Figure~\ref{fig:box_lstsq_experiment} shows the function-value gap across iterations.
    This highlights a commonly observed discrepancy between worst-case convergence guarantees and actual performance in practice.
\end{example}

\subsection{Contributions}
In this paper, we develop a data-driven performance estimation framework for analyzing the convergence of deterministic first-order methods over a distribution of convex optimization problems.
By combining \emph{a priori} information on the function classes and \emph{a posteriori} information from observed algorithm trajectories, our data-driven technique narrows the gap between theoretical and empirical performance.
Our contributions are as follows:
\begin{itemize}
    \item We combine the PEP framework with data-driven Wasserstein distributionally robust optimization (DRO) to create a computational framework that analyzes the performance of a first-order method over a distribution of convex optimization problems.
    To achieve this, we formulate the analysis as a tractable conic program that depends on the observed algorithm trajectories on a finite set of problem instances.
    Our formulation interpolates between in-sample average-case behavior and the worst-case PEP bound as the Wasserstein radius grows, recovering PEP exactly above an explicit threshold.

    \item We establish finite-sample guarantees that our certificate upper-bounds the true expectation and conditional value-at-risk of commonly used performance metrics (\eg, distance to optimality and suboptimality) over unseen instances. 

    \item We show asymptotic rates as well as upper bounds for the conditional value-at-risk of the suboptimality of gradient descent and the fast gradient method on structured distributions of quadratic minimization problems.
    These rates identify settings in which our framework certifies strictly faster performance than the worst case.
    Our rates build on the average-case analyses of~\citet{cunhaOnlyTailsMatter2022,paquetteHaltingTimePredictable2023}, complementing their expectation results with risk-aware guarantees.

    \item Numerical experiments on unconstrained quadratic minimization, logistic regression using a real credit-scoring dataset, and Lasso show that our approach significantly reduces the conservatism of state-of-the-art computer-assisted worst-case bounds.
    Moreover, it directly captures empirical phenomena such as ripples in Nesterov acceleration and heavy-tailed performance variation across problem instances.
    We provide an open-source repository to reproduce all experiments at \url{https://github.com/stellatogrp/dro_pep}.
\end{itemize}

\subsection{Related work}

\paragraph{Classical convergence analysis and algorithm design}
First-order methods have been widely adopted to solve large-scale optimization problems.
Unfortunately, they can require a large number of iterations to reach high-quality solutions, especially for badly scaled problems.
This limitation has led to significant work on the worst-case performance of first-order methods, whose convergence depends critically on the underlying problem structure.
For instance, proximal and fixed-point methods can converge linearly or superlinearly under suitable error-bound conditions~\citep{liaoErrorBoundsPL2024,themelisSuperMannSuperlinearlyConvergent2019}.
Without such structure, convergence is typically sublinear~\citep{bubeckConvexOptimizationAlgorithms2015}, as for the proximal point method~\citep{guTightSublinearConvergence2020,guTightErgodicSublinear2024}, the iterative shrinkage-thresholding algorithm (ISTA), PDHG~\citep{chambolleFirstOrderPrimalDualAlgorithm2011}, and ADMM on nonsmooth convex problems.

To mitigate this slow convergence behavior, Nesterov introduced the fast gradient method (FGM)~\citep{nesterovMethodUnconstrainedConvex1983}, which employs an auxiliary iterate sequence to improve worst-case convergence rates.
Several extensions build on this concept:
the fast iterative shrinkage-thresholding algorithm (FISTA) for composite optimization~\citep{beckFastIterativeShrinkageThresholding2009}, the accelerated proximal point method for monotone inclusion~\citep{kimAcceleratedProximalPoint2021}, the Halpern iteration for fixed-point problems~\citep{liederConvergenceRateHalperniteration2021}, and the extra anchored gradient method for minimax optimization~\citep{yoonAcceleratedAlgorithmsSmooth2021}.
Most of these acceleration schemes are \emph{optimal}, in the sense that their convergence rate matches the fastest achievable rate among first-order methods~\citep{nemirovskyProblemComplexityMethod1983}, with sharper complexity results available for specific setups such as smooth strongly convex minimization~\citep{droriOracleComplexitySmooth2022} and composite minimization with a proximal oracle~\citep{jangComputerAssistedDesignAccelerated2025}.
Nevertheless, worst-case optimality does not always translate to practical efficiency; suboptimal algorithms can outperform accelerated variants on certain problem instances~\citep{parkExactOptimalAccelerated2022}.
To remedy this, we aim to develop a probabilistic performance analysis framework that takes into account the distribution of optimization problems.

\paragraph{Computer-assisted performance analysis}
Several recent approaches cast worst-case analysis as an optimization problem, whose objective is the worst-case convergence rate over a given problem class.
The integral quadratic constraint (IQC) framework uses control theory to find asymptotic convergence rates for strongly convex minimization~\citep{lessardAnalysisDesignOptimization2016,zhangUnifiedAnalysisFirstOrder2021}.
The~\emph{performance estimation problem} uses SDP tools to find finite-step worst-case bounds~\citep{droriPerformanceFirstorderMethods2014,taylorSmoothStronglyConvex2017}, and has been further extended to first-order methods with noisy oracles~\citep{barréPrincipledAnalysesDesign2023,gösgensSubgradientMethodsNonsmooth2025}.
The PEP framework can also be used for optimal step size search, leading to the discovery of worst-case optimal first-order methods including OGM~\citep{kimOptimizedFirstorderMethods2016} and OGM-G~\citep{kimOptimizingEfficiencyFirstOrder2021} for smooth convex minimization, OptISTA~\citep{jangComputerAssistedDesignAccelerated2025} for composite minimization, and APPM for monotone inclusions~\citep{kimAcceleratedProximalPoint2021}.
For parametric linear and quadratic problems, a mixed-integer linear program (MILP) called the verification problem (VP) can evaluate the exact worst-case convergence bounds~\citep{ranjanExactVerificationFirstOrder2025,ranjanVerificationFirstorderMethods2025}.
Although these approaches rely on numerical solvers and yield numerical, non-analytical results,
Lyapunov analysis enables the reconstruction of human-interpretable proofs~\citep{bansalPotentialFunctionProofsGradient2019,upadhyayaAutomatedTightLyapunov2025}.
Still, these computer-assisted techniques provide bounds only for worst-case scenarios.
In this paper, we extend the PEP framework to automatically construct probabilistic performance guarantees of first-order methods.

\paragraph{Probabilistic convergence analysis for deterministic algorithms}
Several works deviate from classical worst-case analysis to study the performance of an algorithm over a distribution of optimization problems.
Notable examples include the average-case analysis~\citep{borgwardtSimplexMethodProbabilistic1987,smaleAverageNumberSteps1983}
and smoothed analysis~\citep{bachOptimalSmoothedAnalysis2025,dadushFriendlySmoothedAnalysis2020,spielmanSmoothedAnalysisAlgorithms2004} of the simplex method, as well as average-case analyses of sorting algorithms~\citep{jiangAveragecaseAnalysisAlgorithms2000,knuthArtComputerProgramming1998,vitanyiAnalysisSortingAlgorithms2007}.
More recent works study the average-case performance of first-order methods, primarily on convex quadratic optimization~\citep{cunhaOnlyTailsMatter2022} and regression problems with Gaussian-distributed features~\citep{paquetteHaltingTimePredictable2023,pedregosaAveragecaseAccelerationSpectral2020,scieurUniversalAverageCaseOptimality2020}.
This literature also extends to accelerated gradient methods for bilinear optimization~\citep{domingo-enrichAveragecaseAccelerationBilinear2021} and distributed consensus algorithms on graphs~\citep{nguyenAveragecaseOptimizationAnalysis2024}.
Beyond average-case, high-probability convergence bounds have been established for PDHG applied to linear programs with sub-Gaussian distributions~\citep{xiongHighProbabilityPolynomialTimeComplexity2025}.
Data-driven approaches in the learning-to-optimize literature~\citep{chenLearningOptimizeTutorial2024} certify an algorithm's performance either through PAC-Bayes generalization bounds~\citep{sambharyaLearningWarmStartFixedPoint2024,sambharyaDataDrivenPerformanceGuarantees2025} or the scenario approach~\citep{huangDataDrivenPerformanceGuarantees2025}.
Most of this literature characterizes average-case or high-probability convergence; performance at specific quantiles has received less attention.
In this paper, we provide a systematic framework that estimates this quantile information through the conditional value-at-risk (CVaR) of a selected performance metric.

\paragraph{Distributionally robust optimization (DRO)}
Distributionally robust optimization is a popular approach to modeling optimization problems affected by random uncertainty, without full knowledge of the underlying distribution.
By constructing an \emph{ambiguity set} of distributions containing the true distribution with high probability, DRO provides a framework to obtain solutions that are robust against distributional uncertainty~\citep{chenDistributionallyRobustLearning2020,kuhnDistributionallyRobustOptimization2025}.
Common ambiguity sets are based on either \emph{a priori} information, such as moments (\eg, Chebyshev, Gelbrich, and Markov ambiguity sets)~\citep{delageDistributionallyRobustOptimization2010,gohDistributionallyRobustOptimization2010,zymlerDistributionallyRobustJoint2013}, or \emph{a posteriori} information, such as sample realizations.
The latter ambiguity sets form balls around the empirical distribution under the Kullback--Leibler divergence or the Wasserstein distance~\citep{mohajerinesfahaniDatadrivenDistributionallyRobust2018,wiesemannDistributionallyRobustConvex2014}.
In this paper, we construct a PEP formulation that uses data-driven Wasserstein DRO~\citep{gaoDistributionallyRobustStochastic2023,kuhnWassersteinDistributionallyRobust2019} to analyze deterministic first-order methods over a distribution of optimization problems.
The random problem instances constitute the uncertainty in the PEP formulation.
Using sampled algorithm trajectories and a support model, we obtain high-confidence probabilistic performance guarantees for unseen instances without imposing a parametric distribution.

\section{Probabilistic analysis of deterministic first-order methods}
\label{sec:prob_analysis_fom}
Consider the convex optimization problem
\begin{equation}\label{prob:base_opt}
    \begin{array}[t]{ll}
        \textnormal{minimize} & f(x),
    \end{array}
\end{equation}
where~$x \in \reals^d$ is the optimization variable and~$f \colon \reals^d \to \reals$ is a convex objective function.
For simplicity, we focus here on the smooth minimization case, \ie, the case where $f$ is differentiable with $L$-Lipschitz continuous gradient for some~$L > 0$.
We assume that this problem has an optimal solution~$x^\star \in \reals^d$ with optimal value $f^\star = f(x^\star)$ satisfying the stationarity condition~$g^\star = \nabla f(x^\star) = 0$.

In \suppref{appendix:composite}, we extend our analysis to composite minimization problems with a nonsmooth, possibly extended-valued term, solved via proximal algorithms.

\paragraph{Algorithm}
Given an initial iterate $x^0\in\reals^d$, we refer to the mapping~$\algo \colon (f, x^0) \mapsto \{x^k\}_{k=0,1,\dots}$ as an \emph{algorithm}.
We say that~$\algo$ is a \emph{fixed-step first-order method} if there exists a sequence of scalar step sizes~$\{\eta_k^i\}_{i=0}^{k}$ for $k=0,1,\dots$ such that
\begin{equation*}
    x^{k+1} = x^k - {\textstyle \sum_{i=0}^{k}} \eta_k^i\, g^i,
    \qquad k=0,1,\dots,
\end{equation*}
where~$g^i = \nabla f(x^i)$ is the gradient of $f$ evaluated at $x^i$.

\paragraph{Performance metric}
We measure the performance of the algorithm~$\algo$ applied to~$f$ from the initial iterate~$x^0$ using a scalar-valued function~$\perfmetr[K](f, x^0)$ that quantifies the behavior of the iterates~$\{x^k\}_{k=0}^K$ up to iteration~$K$.
Examples of performance metrics include the squared distance to optimality~$\|x^K - x^\star\|^2$, the squared gradient norm~$\|\nabla f(x^K)\|^2$, and the function-value suboptimality~$f(x^K) - f^\star$.
For the finite-step performance estimation framework developed in Section~\ref{sec:convex_formulation}, we assume a~\emph{large-scale setting} where the underlying problem dimension~$d$ is substantially larger than~$K$.

\paragraph{Function class}
We refer to a collection~$\fclass$ of functions $f$ with shared characteristics as a \emph{function class}.
For~$0 \le \mu \le L$, the function classes we consider are:
\begin{itemize}
\item $\mathcal{F}_{0,L}$: $L$-smooth convex functions.
    \item $\mathcal{F}_{\mu,L}$: $L$-smooth $\mu$-strongly convex functions.
    \item $\mathcal{Q}_{\mu,L}$: convex quadratic functions $x \mapsto (1/2) x^T Q x$ with $\mu I \preceq Q \preceq L I$.
\end{itemize}

\paragraph{Initial condition}
To meaningfully bound the performance metric $\perfmetr[K](f, x^0)$, we restrict the initial iterate~$x^0$ to be in the set~$\xclass(f) = \{ x \in \reals^d \mid \|x - x^\star\|^2 \le r^2 \} \subseteq \reals^d$, where the distance to an optimal solution is bounded by $r>0$.

We refer to the set of all such function and initial iterate pairs~$(f, x^0)$ as~$\fclassX = \{ (f, x^0) \in \fclass \times \reals^d \mid x^0 \in \xclass(f) \}$.

\subsection{Risk measures for performance analysis}\label{subsec:risk_measure}
Let~$\measure(\fclassX)$ be the set of all probability distributions supported on~$\fclassX$.
We model the distribution of the optimization problem instances and the initial iterates as~${(f, x^0) \sim \prob}$, where~$\prob \in \measure(\fclassX)$.
To derive convergence guarantees with respect to~$\prob$, we consider risk measures of~$\perfmetr[K](f, x^0)$.
The \emph{value-at-risk} ($\var$) with quantile parameter~$\alpha \in (0, 1)$ is defined as~$\var_{\alpha} (\perfmetr[K](f, x^0)) = \inf \big\{\tau \in \reals \mid \prob(\perfmetr[K](f, x^0) \le \tau) \ge 1 - \alpha \big\}$,
and is exactly the~$(1-\alpha)$-quantile of $\perfmetr[K](f, x^0)$ over $(f, x^0) \sim \prob$.
Direct optimization involving $\var$ is generally nonconvex~\citep{uryasevConditionalValueatRiskOptimization2001}.
A tractable convex surrogate is the \emph{conditional value-at-risk} ($\cvar$), defined for~$\alpha \in (0, 1]$ as
\begin{equation*}
    \cvar_\alpha(\perfmetr[K](f, x^0))
    =
    \inf_{t \in \reals}\,  \Big(
        t + \alpha^{-1} \Expect_{(f, x^0) \sim \prob} (\perfmetr[K](f, x^0) - t)_+
    \Big),
\end{equation*}
where~$(z)_+ = \max \{z, 0\}$.
By~\citet[Theorem~6.2]{shapiroChapter6Risk2021}, for~$\alpha \in (0, 1)$, $\cvar_\alpha$ is the average of the upper~$\alpha$-tail of~$\perfmetr[K](f, x^0)$, with partial mass at~$\var_\alpha$ when necessary; it consequently upper-bounds $\var_\alpha$.

By varying the quantile parameter~$\alpha$,
we can interpolate between the average case ($\alpha = 1$), \ie, $\cvar_{1}(\perfmetr[K](f, x^0)) = \Expect_{(f, x^0) \sim \prob} (\perfmetr[K](f, x^0))$,
and the worst case ($\alpha \to 0^+$), \ie, $\lim_{\alpha \to 0^+} \cvar_\alpha(\perfmetr[K](f, x^0)) = \esssup_{(f, x^0)\in\fclassX}\,\perfmetr[K](f, x^0)$.
Here the essential supremum is taken with respect to~$\prob$, and equals $\sup_{(f, x^0) \in \fclassX} \perfmetr[K](f, x^0)$ when~$\prob$ has full support on~$\fclassX$~\citep{uryasevConditionalValueatRiskOptimization2001}.
We measure performance using the risk measures summarized in Table~\ref{tab:perfmeas}.
\begin{table}[ht]
    \centering
    \caption{
        List of performance measures with respect to $(f, x^0) \sim \prob$.
All measures can be expressed as a function of $\cvar_\alpha(\perfmetr[K](f, x^0))$ with specific values of~$\alpha \in (0, 1]$.
    }
    \begin{tabular}{lll}
        \toprule
            &   performance measure  & quantile parameter~$\alpha$ \\
        \midrule 
            average case &
                $   
                    \Expect_{(f, x^0) \sim \prob} \left(\perfmetr[K](f, x^0)\right)
                $ & $\alpha = 1$ \\
            $\cvar$ &
                $   
                    \cvar_\alpha(\perfmetr[K](f, x^0))
                $ & $\alpha \in (0, 1)$ \\
            worst case  &
                $   
                    \esssup_{(f, x^0)\in\fclassX}\, \perfmetr[K](f, x^0)
                $ & $\alpha \to 0^+$ \\
        \bottomrule
    \end{tabular}
    \label{tab:perfmeas}
\end{table}

Evaluating these risk measures requires solving the minimization problem
\begin{equation}\label{eq:expect_form}
    \inf_{\ell \in \lclass}\, \Expect_{(f, x^0) \sim \prob}\, \left( \ell(f, x^0) \right),
\end{equation}
over the \emph{loss function}~$\ell \colon \fclassX \to\reals$.
Specifically, we set~$\lclass = \{(f, x^0) \mapsto \perfmetr[K](f, x^0)\}$ for the average case and~$\lclass = \{(f, x^0) \mapsto t + \alpha^{-1}(\perfmetr[K](f, x^0) - t)_+ \mid t \in \reals \}$ for the~$\cvar$ case.
Unfortunately, solving~\eqref{eq:expect_form} presents two main challenges.
First, it requires representing functions living in the infinite-dimensional space~$\fclass$.
Second, evaluating~$\Expect_{(f, x^0)\sim\prob}$ requires full knowledge of~$\prob$, which is rarely available in practice.
Section~\ref{sec:convex_formulation} addresses both issues.

\subsection{Probabilistic convergence rates for quadratic minimization}
\label{subsec:prob_rates}
Before developing our data-driven formulation, we study a setting where exact asymptotic convergence rates are known and probabilistic analysis provably leads to faster convergence than the worst case.
We consider the minimization problem~\eqref{prob:base_opt} with quadratic objective
\begin{equation}\label{eq:quad_obj}
    f(x) = (1/2) (x - x^\star)^T Q (x - x^\star),
\end{equation}
where~$Q \in \symm_+^d$ has eigenvalues in~$[0, L]$ and~$x^\star$ is an optimal solution.
We draw random instances of~\eqref{eq:quad_obj} in fixed dimension~$d$ as follows.
This asymptotic subsection uses a separate, possibly unbounded instance model and does not impose the uniform initial-radius condition from the finite-step framework.
The eigenvalues of~$Q$ are exchangeable (\eg, i.i.d.), so each follows a common marginal distribution~$\nu$, whose density behaves as~$p(\lambda) \sim p_0 \lambda^a$ near~$\lambda = 0$ for some~$a > -1$ and~$p_0>0$.
The initial error~$x^0 - x^\star$ is independent of~$Q$ and rotationally invariant, \ie, it has no preferred direction, with~$\Expect \| x^0 - x^\star \|^2 = r^2$.

For this instance distribution, the average-case analyses of~\citet{scieurUniversalAverageCaseOptimality2020,cunhaOnlyTailsMatter2022,paquetteHaltingTimePredictable2023} derive exact asymptotic rates for the expected performance of the following first-order methods:
gradient descent (GD)
\begin{equation}\label{eq:GD}\tag{GD}
    x^{k+1} = x^{k} - \eta \nabla f(x^k),
    \qquad k=0,1,\dots,
\end{equation}
and Nesterov's fast gradient method (FGM)~\citep{nesterovMethodUnconstrainedConvex1983},
\begin{equation}\label{eq:FGM}\tag{FGM}
    \begin{array}{ll}
        x^{k+1} &= y^k - \eta \nabla f(y^k)\\
        y^{k+1} &= x^{k+1} + \gamma_k (x^{k+1} - x^k),
    \end{array}
    \qquad k=0,1,\dots,
\end{equation}
where~$y^0 = x^0$, $\eta \in (0, 1/L]$, and~$\gamma_k \ge 0$ is a momentum coefficient; in this section, we take~$\gamma_k = k / (k+3)$.
Their average-case performances are faster than the respective worst-case rates: $\Theta(1/K)$ for GD~\citep{kimProofExactConvergence2025} and $\Theta(1/K^2)$ for FGM~\citep{daspremontAccelerationMethods2021}.
Building on these results, the following proposition complements the average-case performance with the~$\cvar$-performance introduced in Section~\ref{subsec:risk_measure}.
In contrast to the large-scale finite-step regime considered in the rest of the paper, this result fixes the dimension~$d$ and lets~$K \to \infty$.
The proof is in Appendix~\ref{appendix:ls_rate}.

\begin{proposition}[Probabilistic convergence rate of quadratic minimization]
    \label{prop:ls_rate}
    Consider random instances of problem~\eqref{prob:base_opt} with the quadratic objective~\eqref{eq:quad_obj}, generated as described above.
    Fix~$\alpha \in (0, 1]$ and let the performance metric be~$\perfmetr[K] = f(x^K) - f^\star$, with GD and FGM run at step size~$\eta = 1/L$.
    As~$K \to \infty$, the following statements hold:
    \begin{enumerate}
        \item[(i)]
        The $\cvar$-performance is asymptotically equivalent to the scaled average-case performance:
        \begin{equation*}
            \cvar_\alpha(\perfmetr[K]) \sim \alpha^{-1} \Expect(\perfmetr[K]),
        \end{equation*}
        where~$\Expect(\perfmetr[K]) = \Theta\big(K^{-(a+2)}\big)$ for GD and~$\Expect(\perfmetr[K]) = \Theta\big(K^{-e(a)}\big)$ for FGM, with~$e(a) = 2(a+2)$ for~$a \le -1/2$ and $e(a) = a + 7/2$ for~$a \ge -1/2$ with an additional~$\log K$ factor at~$a = -1/2$.

        \item[(ii)]
        Suppose that~$(x^0 - x^\star)$ is furthermore Gaussian and~$\Expect (w_i w_j) \le \Expect (w_i) \Expect (w_j)$ for~$i \ne j$,
        where~$w_i = \lambda_i \rho_K (\lambda_i)^2$ and the residual polynomial~$\rho_K$ is defined by~$x^K-x^\star=\rho_K(Q)(x^0-x^\star)$ for GD or FGM\@.
        Then the~$\cvar$-performance satisfies
        \begin{equation*}
            \cvar_\alpha(\perfmetr[K])
            \le
                \sqrt{\alpha^{-1} \Expect\big((\perfmetr[K])^2\big)}
                = O \left( \alpha^{-1/2}\left(d^{-1/2} K^{-h(a)} + \Expect(\perfmetr[K])\right) \right),
        \end{equation*}
        where~$h(a) = (a+3)/2$ for GD, and for FGM, $h(a) = a+3$ if $a \le 0$ and $h(a) = a/2 + 3$ if $a \ge 0$, with an additional~$\sqrt{\log K}$ multiplier at~$a = 0$.
        For each fixed~$d$, the expectation term is lower order, so the bound simplifies to~$O \big( \alpha^{-1/2} d^{-1/2} K^{-h(a)} \big)$; this simplification is not uniform in~$d$.
    \end{enumerate}
\end{proposition}

\begin{remark}[On the cross-moment condition]
    The condition~$\Expect (w_i w_j) \le \Expect (w_i) \Expect (w_j)$ for~$i \ne j$ in part~($ii$) requires the terms~$w_i$ to be pairwise non-positively correlated.
    For every fixed~$K$, it holds asymptotically with equality as~$d\to\infty$ for a uniformly sampled pair of distinct eigenvalues in the Marchenko--Pastur limit, as in the Wishart model of Section~\ref{subsec:quad_experiment}.
    It holds exactly with equality when the eigenvalues are sampled independently from~$\nu$.
\end{remark}

Proposition~\ref{prop:ls_rate} also illustrates why estimating probabilistic performance from a small sample can be challenging.
Given an integer~$K\ge1$, the spectral loss kernel~$\lambda \mapsto \lambda (1-\lambda / L)^{2K}$ of GD peaks at~$\lambda_\star = L/(2K+1)$, so problem instances with an eigenvalue at scale~$\Theta(1/K)$ contribute most to the mean.
However, a small sample may not include such instances, making the empirical estimate of~$\perfmetr[K]$ overly optimistic.
This sampling difficulty motivates the Wasserstein DRO formulation developed next, which hedges against undersampled instances while providing a high-confidence bound.

\section{Convex formulation of the data-driven performance estimation problem}
\label{sec:convex_formulation}
In this section, we derive a tractable convex formulation that upper-bounds the expected loss~\eqref{eq:expect_form} with high probability.
We address the two challenges raised in Section~\ref{sec:prob_analysis_fom} by replacing the infinite-dimensional function class~$\fclass$ by finite-dimensional interpolation conditions,
and the unknown distribution~$\prob$ by an empirical distribution with a Wasserstein ambiguity set.

\subsection{Finite-dimensional representation via interpolation conditions}
\label{subsec:finite_formulation}
We obtain finite-dimensional representations of~$(f, x^0) \in \fclassX$ by evaluating the function~$f$ and its gradient only at the points encountered by the algorithm starting from~$x^0$, which is the only function information required to compute the performance metrics.

Suppose we are given a sequence of vectors~$z^k = (x^k, g^k, f^k) \in \reals^{2d+1}$ for~$k\in \mathcal{I}$,
where $\mathcal{I}$ is an arbitrary index set.
We say that~$\{z^k\}_{k \in \mathcal{I}}$ is \emph{$\fclass$-interpolable} if there exists a function~$f \in \fclass$ such that~$f^k = f(x^k)$ and~$g^k = \nabla f(x^k)$ for all~$k \in \mathcal{I}$.
For each function class in Section~\ref{sec:prob_analysis_fom}, interpolability is characterized by necessary and sufficient quadratic inequalities in the entries of~$\{z^k\}_{k\in \mathcal{I}}$. These inequalities are called the~\emph{interpolation conditions}~\citep{droriPerformanceFirstorderMethods2014,ryuOperatorSplittingPerformance2020,taylorExactWorstCasePerformance2017,taylorSmoothStronglyConvex2017}.

In this paper, we consider the tuples at the algorithm iterates and at the optimal solution~$x^\star$, which correspond to the index set~$\mathcal{I}=\{\star, 0, 1, \dots, K\}$.
Note that, without loss of generality, we can set~$x^\star = 0$ and~$f^\star = 0$, due to translation invariance of the function classes considered~\citep{taylorSmoothStronglyConvex2017}.

\paragraph{Semidefinite variable lifting}
Being quadratic in the entries of the tuples~$\{z^k\}_{k \in \mathcal{I}}$, the interpolation conditions are nonconvex in general.
Adopting the strategy used in the PEP framework~\citep{droriPerformanceFirstorderMethods2014,taylorExactWorstCasePerformance2017,taylorSmoothStronglyConvex2017}, we lift the algorithm iterates to a higher-dimensional space where they can be represented as a set of convex constraints.

Given a set of tuples~$\{z^k\}_{k \in \mathcal{I}}$ with~$z^k = (x^k, g^k, f^k)$, define the symmetric positive semidefinite \emph{Gram matrix}~$G \in \symm_+^{K+2}$ and the function-value vector~$F \in \reals^{K+1}$ as
\begin{equation}
        G =
        P^T P,
            \qquad
        F =
        \left(
            f^0 - f^\star, \dots, f^K - f^\star
        \right),
    \label{eq:gram-matrix}
\end{equation}
where~$P = \begin{bmatrix} (x^0 - x^\star) & g^0 & \cdots & g^K \end{bmatrix} \in \reals^{d \times (K+2)}$.
The interpolation conditions for the function classes in Section~\ref{sec:prob_analysis_fom} can be written as linear inequalities in~$G$ and~$F$~\citep{taylorSmoothStronglyConvex2017}.
For simplicity, we write such inequalities,
as well as the performance metric and the initial condition,
in terms of the inner product~$\inner{\cdot}{\cdot}$ defined as~$\langle(X, Y),\, (G, F) \rangle = \Tr(X^T G) + Y^T F$ for any~$(X, Y),\, (G, F) \in \symm^{K+2} \times \reals^{K+1}$.

\begin{assumption}
    \label{assumption:interpolation}
    The interpolation conditions of the function class~$\fclass$ are a set of linear matrix inequalities of the form
    \begin{equation*}
        \inner{(A_m, b_m)}{(G, F)} \le 0,
        \qquad m \in M = \{ (p, q) \in \mathcal{I} \times \mathcal{I} \mid p \neq q \},
    \end{equation*}
    where~$A_m \in \symm^{K+2}$ and~$b_m \in \reals^{K+1}$.
    By an abuse of notation, we reuse~$\fclass$ and~$\xclass$ for their lifted representations. Accordingly, we write
    \begin{equation*}
        \fclass = \left\{
            (G, F) \in \symm_+^{K+2} \times \reals^{K+1}
            \;\middle|\;
                \linear(G, F) \in \reals_+^{|M|}
        \right\},
    \end{equation*}
in terms of the linear operator~$\linear \colon \symm^{K+2} \times \reals^{K+1} \to \reals^{|M|}$ and its adjoint~$\linear^*$ defined as
    \begin{equation*}
        \linear(G, F) 
        =
        \big(
            - \inner{(A_m, b_m)}{(G, F)}
        \big)_{m \in M} \,,
        \qquad
        \linear^*(y)
        = - {\textstyle \sum_{m \in M}} y_m (A_m, b_m).
    \end{equation*}
    The performance metric~$\perfmetr[K](G, F)$ has the form~$\perfmetr[K](G, F) = \inner{(\Aobj, \bobj)}{(G, F)}$,
    where~$\Aobj \in \symm^{K+2}$ and~$\bobj \in \reals^{K+1}$.
    Lastly, the initial condition~$x^0 \in \xclass(f)$ for~$f \in \fclass$ has the form~$\inner{(A_0, b_0)}{(G, F)} + c_0 \le 0$, where~$A_0 \in \symm^{K+2}$, $b_0 \in \reals^{K+1}$, and~$c_0 \in \reals$.
    We may omit the dependence on~$f$ and write it as
    \begin{equation*}
        \xclass =
            \left\{
                (G, F) \in \symm_+^{K+2} \times \reals^{K+1}
                \;\middle|\;
                    \inner{(A_0, b_0)}{(G, F)} + c_0 \le 0
            \right\},
    \end{equation*}
    since the definition of $G = P^T P$ in~\eqref{eq:gram-matrix} centers the first column of $P$ around $x^\star$.
\end{assumption}

By construction, the transformation~$(f, x^0) \mapsto (G, F)$ maps~$\fclassX$ into~$\fclass \cap \xclass$, since the tuples~$\{z^k\}_{k \in \mathcal{I}}$ are generated by~$f \in \fclass$ and the initial iterate satisfies~$x^0 \in \xclass(f)$.
Note that through this transformation, we lose information about the dimension~$d$ of the underlying problem.
Nevertheless, as long as~$d$ is large enough, \ie, $d \ge K+2$, the finite-dimensional representation~$(G, F)$ is enough to express the behavior of the given algorithm on~$(f, x^0) \in \fclassX$~\citep[Theorems~5 and~6]{taylorSmoothStronglyConvex2017},
which also underlies the strong duality argument crucial for the tightness of the worst-case analysis in the PEP framework.

\paragraph{From function distances to Gram-matrix distances}
We now relate the distance between problem instances~$(f, x^0)$ to that between their lifted representations~$Z = (G, F)$, focusing on the function class~$\fclass$ with Lipschitz gradients.
Since the iterates of~$\algo$ over~$K$ steps remain in the convex compact set~$\domain$ from Lemma~\ref{lem:domain_size} (deferred to \suppref{appendix:trajectory_implication}),
the distance between two functions~$f_1, f_2 \in \fclass$ defined in~\eqref{eq:c1_sup_norm} is well-defined and finite:
\begin{equation}\label{eq:c1_sup_norm}
    \norm{f_1 - f_2}
    =
        \sup_{x \in \domain}\, |f_1(x) - f_2(x)|
        + \sup_{x \in \domain}\, \norm{\nabla f_1(x) - \nabla f_2(x)}.
\end{equation}
The expression in~\eqref{eq:c1_sup_norm} is a metric on restrictions of functions to~$\domain$, or equivalently on functions identified whenever they agree on~$\domain$.
The distance between~$Z_1 = (G_1, F_1)$ and~$Z_2 = (G_2, F_2)$ in~$\fclass \cap \xclass$ is defined as~$\|Z_1 - Z_2\| = \sqrt{ \|G_1 - G_2\|_F^2 + \|F_1 - F_2\|_2^2 }$, with Frobenius norm~$\|\cdot\|_F$.
We denote its dual norm by~$\|\cdot\|_*$ that is a self-dual Euclidean product norm.

We claim that whenever two function--initial-iterate pairs are close, their lifted representations are also close.
See \suppref{appendix:trajectory_implication} for the proof of the general composite optimization setting, which reduces to the single-objective case~\eqref{prob:base_opt} below.
\begin{theorem}\label{thm:trajectory_implication}
Suppose that~$\fclass$ is either~$\mathcal{F}_{0, L}$, $\mathcal{F}_{\mu, L}$, or $\mathcal{Q}_{\mu, L}$.
Let~$\algo$ be a fixed-step first-order method with iteration budget~$K \ge 1$.
    There exists~$C_K > 0$ such that
    \begin{equation*}
        \norm{Z_1 - Z_2}
        \le
            C_K \| (f_1, x_1^0) - (f_2, x_2^0) \|
            = C_K \big(
                \|f_1 - f_2\|
                + \|x_1^0 - x_2^0\|
            \big),
    \end{equation*}
    for any~$x_1^0, x_2^0 \in \domain$ and continuously differentiable~$f_1, f_2 \in \fclass$, under the normalization~$x_1^\star=x_2^\star=0$ and~$f_1^\star=f_2^\star=0$.
\end{theorem}

\subsection{Approximating the expected loss from data}
\label{subsec:DRO}
After the variable lifting from Section~\ref{subsec:finite_formulation}, the expectation in~\eqref{eq:expect_form} equals
\begin{equation}\label{eq:true_exp_z}
    \Expect_{Z \sim \prob_Z} \left( \ell (Z) \right),
\end{equation}
where~$Z = (G, F)$ and~$\prob_Z \in \measure(\fclass \cap \xclass)$ is the pushforward of~$\prob \in \measure(\fclassX)$ under the mapping~$(f, x^0) \mapsto Z$.
Unfortunately, in most practical scenarios, complete knowledge of~$\prob_Z$ is unavailable.
We address this issue using \emph{distributionally robust optimization}~\citep{kuhnWassersteinDistributionallyRobust2019,kuhnDistributionallyRobustOptimization2025,mohajerinesfahaniDatadrivenDistributionallyRobust2018}.
Specifically, we incorporate \emph{a priori} structural information modeled in the previous sections, such as initial conditions and the interpolation conditions, with data-driven \emph{a posteriori} information from running the algorithm on sample realizations of $(f, x^0) \sim \prob$.

\paragraph{Distributionally robust optimization}
Consider a set~$\mathcal{P}$ of probability distributions supported on~$\fclass \cap \xclass$.
If~$\prob_Z \in \mathcal{P}$, then the worst-case expectation over~$\mathcal{P}$ upper-bounds the true expectation~\eqref{eq:true_exp_z} over~$\prob_Z$, \ie,
\begin{equation}\label{eq:worst-case-exp}
    \inf_{\ell \in \lclass}\, \Expect_{(f, x^0) \sim \prob} \left( \ell(f, x^0) \right)
    =
    \inf_{\ell \in \lclass}\,\underset{Z \sim \prob_Z}{\Expect} \left( \ell(Z) \right)
    \le
    \begin{array}[t]{ll}
        \underset{\probQ_Z \in \mathcal{P}}{\sup}
        &    \underset{\ell \in \lclass}{\inf}\, \underset{Z \sim \probQ_Z}{\Expect} \left( \ell(Z) \right).
    \end{array}
\end{equation}
Our goal is to construct an \emph{ambiguity set}~$\mathcal{P}$ that is guaranteed to contain~$\prob_Z$ with high probability,
using only partial information about~$\prob_Z$ obtained from samples~\mbox{$\hZ_1,\dots,\hZ_N \sim \prob_Z$}.

\paragraph{Finite-sample guarantee of Wasserstein DRO}
Consider the empirical distribution~$\hprob_Z = (1/N) \sum_{i=1}^N \delta_{\hZ_i}$, where~$\delta_{\hZ_i}$ is the Dirac distribution at the~$i$-th sample~$\hZ_i = (\hG_i, \hF_i)$.
We set the ambiguity set~$\mathcal{P}$ as the ball of radius~$\varepsilon$ around~$\hprob_Z$:
\begin{equation*}
    \ambiset[\varepsilon]
    = \left\{
        \probQ_Z \in \measure(\fclass \cap \xclass)
        \;\middle|\;
            W \left( \hprob_Z, \probQ_Z \right) \le \varepsilon
    \right\},
\end{equation*}
with respect to the $1$-Wasserstein distance
\begin{align*}
    W (\prob_Z, \probQ_Z)
    &=
        \inf_{\pi \in \Pi(\prob_Z, \probQ_Z)}\,
        \biggl(
            \int_{(\fclass \cap \xclass) \times (\fclass \cap \xclass)} \norm{Z_1 - Z_2}
            \,\pi(dZ_1, dZ_2)
        \biggr),
\end{align*}
where~$\pi$ is a joint distribution of $Z_1$ and~$Z_2$ with marginals $\prob_Z$ and $\probQ_Z$.

Given a failure probability~$\beta \in (0, 1)$,
well-known measure concentration results~\citep[Theorem~2]{fournierRateConvergenceWasserstein2015} and~\citep[Proposition~20]{weedSharpAsymptoticFinite2019} guarantee the existence of a radius~$\varepsilon = \varepsilon_N (\beta) >0$ such that the true distribution $\prob_Z$ lies in~$\ambiset[\varepsilon]$ with probability at least~$1-\beta$, \ie,
\begin{equation}\label{eq:high_prob_guarantee}
\prob_Z^N \left( W ( \prob_Z, \hprob_Z ) \le \varepsilon \right) \ge 1 - \beta.
\end{equation}
Consequently, the worst-case expectation in~\eqref{eq:worst-case-exp} with~$\mathcal{P} = \ambiset[\varepsilon]$ upper-bounds the true expectation~\eqref{eq:true_exp_z} with probability at least~$1-\beta$.
Note that the support set~$\fclass \cap \xclass$ of the distributions considered is compact,
as the algorithm trajectories are uniformly bounded over~$(f, x^0) \in \fclassX$ from Lemma~\ref{lem:domain_size},
given $x^\star = 0$, $f^\star = 0$, and $\|x^0 - x^\star\| \le r$ uniformly over~$(f, x^0) \in \fclassX$.
This compactness permits the finite-sample guarantee~\eqref{eq:high_prob_guarantee}, because the support set is bounded and has finite covering numbers.

\paragraph{Connection to worst-case expectation over infinite-dimensional function spaces}
We now interpret the worst-case value in the right-hand side of~\eqref{eq:worst-case-exp} as an upper bound on the worst-case expectation over the infinite-dimensional space~$\fclassX$ of functions and initial iterates.
The proof in Appendix~\ref{appendix:worst-case-expectations} combines the Wasserstein convergence of empirical measures on compact metric spaces with Theorem~\ref{thm:trajectory_implication}.
Let~$W_{\fclassX}$ denote the~$1$-Wasserstein distance induced by the instance metric~$\|(f_1,x_1^0)-(f_2,x_2^0)\|$ in Theorem~\ref{thm:trajectory_implication}.

\begin{proposition}
    \label{prop:worst-case-expectations}
Consider the setup of Theorem~\ref{thm:trajectory_implication} with~$C_K > 0$.
    Given a failure probability~$\beta \in (0, 1)$,
    there exists~$\varepsilon = \varepsilon_N (\beta) >0$ such that the instance-space Wasserstein ball~$\mathcal{P}=\mathcal{B}_{\fclassX}(\hprob,\varepsilon/C_K):=\{\,\probQ\in\measure(\fclassX)\mid W_{\fclassX}(\probQ,\hprob)\le\varepsilon/C_K\,\}$, centered at~$\hprob = (1/N) \sum_{i=1}^N \delta_{(\hf_i, \hx_i^0)}$, satisfies the following with probability at least~$1-\beta$:
\begin{equation*}
        \Expect_{(f, x^0) \sim \prob} \big( \ell( f, x^0 ) \big)
        \le
        \underset{\probQ \in \mathcal{P}}{\sup}\, \underset{(f, x^0) \sim \probQ}{\Expect} \big( \ell( f, x^0 ) \big)
        \le
        \underset{\probQ_Z \in \ambiset[\varepsilon]}{\sup}\, \underset{Z \sim \probQ_Z}{\Expect} \left( \ell(Z) \right).
    \end{equation*}
\end{proposition}

\subsection{Tractable convex formulation using duality theory}
\label{subsec:tractable}
The following theorem gives a tractable convex formulation of the worst-case expectation on the right-hand side of~\eqref{eq:worst-case-exp},
whose proof is deferred to Section~\ref{subsec:proof-conic-dro-pep}.
\begin{theorem}
    \label{thm:conic-dro-pep}
    Consider the fixed-step first-order method~$\algo$ with integer iteration budget~$K \ge 1$ and its samples~$\{(\hG_i, \hF_i)\}_{i=1}^N$ drawn from~$\prob_Z$.
    Suppose that~$\fclass = \mathcal{F}_{\mu,L}$ for some~$0 \le \mu < L < \infty$ and that the method uses the most recent gradient at every iteration, \ie, $\eta_k^k \neq 0$ for~$k = 0,\dots,K-1$.
    Let the set~$\lclass$ of loss functions~$\ell \colon \fclass \cap \xclass \to \reals$ parametrized by~$t \in \reals$ be
    \begin{equation*}
        \lclass = \left\{
            (G, F) \mapsto \max_{1\le j\le J}\, \ell^j (G, F)
            \;\middle|
                \begin{array}{l}
                    \ell^j (G, F) = \inner{(\Aobj^j, \bobj^j)}{(G, F)} + \cobj^j t, \\
                    j = 1, \dots, J,
                    \quad t \in \reals
                \end{array}
        \right\},
    \end{equation*}
    with~$\Aobj^j \in \symm^{K+2}$, $\bobj^j \in \reals^{K+1}$, and~$\cobj^j \in \reals$.
    Under Assumption~\ref{assumption:interpolation},
    \begin{equation*}
        \underset{\probQ_Z \in \ambiset[\varepsilon]}{\sup}\,
        \underset{\ell \in \lclass}{\inf}\,
        \underset{(G, F) \sim \probQ_Z}{\Expect}\, \left( \ell (G, F) \right) =
        \setlength{\arraycolsep}{2pt}
        \begin{array}[t]{ll}
            \textnormal{maximize}
            &   \sum_{i=1}^N \sum_{j=1}^J \pi_i^j \big\langle (\Aobj^j, \bobj^j),\, (G_i^j, F_i^j) \big\rangle \\
            \textnormal{subject to}
                &   \sum_{i=1}^N \sum_{j=1}^J \pi_i^j \|(G_i^j, F_i^j) - (\hG_i, \hF_i)\| \le \varepsilon,
        \end{array}
    \end{equation*}
    for every~$\varepsilon>0$,
    with variables~$(\pi_i^j) \in \Gamma$ and~$(G_i^j, F_i^j) \in \fclass \cap \xclass$ where
    \begin{align*}
        \Gamma = \left\{
            (\pi_i^j) \in \reals_+^{N J}
            \;\middle|
                \begin{array}{ll}
                    \sum_{j=1}^J \pi_i^j = 1/N,
                        & i=1,\dots,N  \\
                    \sum_{j=1}^J \left( \sum_{i=1}^N \pi_i^j \right) \cobj^j = 0,\\
                \end{array}
        \right\}.
    \end{align*}
    This is also equivalent to the following minimization problem, denoted by~\eqref{prob:dro-pep}:
    \begin{equation*}
        \setlength{\arraycolsep}{2pt}
        \begin{array}[t]{ll}
            \textnormal{minimize}
            &   (1/N) \sum_{i=1}^N s_i \\
            \textnormal{subject to}
            &
                \cobj^j t
                - c_0 \tau_i^j
                - \big\langle (X_i^j, Y_i^j),\, (\hG_i, \hF_i) \big\rangle
                + \lambda \varepsilon \le s_i \\
            &   - \linear^*(y_i^j) - (X_i^j, Y_i^j) + \tau_i^j (A_0, b_0) - (\Aobj^j, \bobj^j) \in \symm_+^{K+2} \times \{0\} \\
            &   \norm[*]{(X_i^j, Y_i^j)} \le \lambda,
                \hfill  i = 1,\dots,N,\;
                        j = 1,\dots,J,
        \end{array}
        \makeatletter\def\@currentlabel{DRO-PEP}\makeatother
        \label{prob:dro-pep}
\end{equation*}
    with variables~$s_i \in \reals$, $t \in \reals$, $\tau_i^j \in \reals_+$, $(X_i^j, Y_i^j) \in \symm^{K+2} \times \reals^{K+1}$, $\lambda \in \reals_+$, and~$y_i^j \in \reals_+^{|M|}$.
\end{theorem}

\begin{remark}[Problem complexity]
    The maximization problem of Theorem~\ref{thm:conic-dro-pep} is a finite-dimensional nonconvex program,
    bilinear in the weights~$(\pi_i^j)$ and the lifted variables~$(G_i^j, F_i^j)$,
    whereas the minimization problem~\eqref{prob:dro-pep} is a convex conic program.
    The problem size of~\eqref{prob:dro-pep} grows in both the number of iterations~$K$ and the number of samples~$N$:
    it has~$NJ$ variables~$(X_i^j, Y_i^j) \in \symm^{K+2} \times \reals^{K+1}$, each appearing in a semidefinite constraint of size~$K+2$.
    For a fixed ambiguity radius~$\varepsilon$, the concentration results behind~\eqref{eq:high_prob_guarantee} permit a smaller certified failure probability~$\beta$ as~$N$ grows.
    Thus, more samples provide higher confidence that the bound~\eqref{eq:worst-case-exp} holds, at the price of a larger problem.
\end{remark}

\begin{remark}[Constructing a worst-case distribution]
    By~\citet[Corollary~4.6]{mohajerinesfahaniDatadrivenDistributionallyRobust2018}, a worst-case distribution achieving the optimal value of~\eqref{prob:dro-pep} is a weighted sum of Dirac distributions.
    Under the standing assumption~$d\ge K+2$, let~$(\pi_i^j) \in \Gamma$ and~$\{(G_i^j, F_i^j)\}_{i,j}$ be optimal solutions to the maximization problem of Theorem~\ref{thm:conic-dro-pep}.
    For each~$(G_i^j, F_i^j)$, there exists an interpolating function~$f_i^j \in \fclass$ and an initial point~$(x_i^j)^0 \in \reals^d$ constructed from the factorization~$G_i^j = (P_i^j)^T P_i^j$~\citep[Remark~2]{taylorSmoothStronglyConvex2017}.
    Note that there can be multiple such functions and initial points.
    Then we can construct a worst-case distribution~$\probQ^\star = \sum_{i=1}^N \sum_{j=1}^J \pi_i^j \delta_i^j$
    where each~$\delta_i^j$ is a Dirac distribution on~$( f_i^j, (x_i^j)^0 ) \in \fclassX$.
\end{remark}

In the following corollary, we specialize Theorem~\ref{thm:conic-dro-pep} to upper-bound the true performance in terms of expectation and $\cvar$.
\begin{corollary}\label{cor:conic-dro-pep}
    Let the failure probability~$\beta\in(0,1)$ and the ambiguity set radius~$\varepsilon = \varepsilon_N(\beta)$ be given as in the finite-sample guarantee~\eqref{eq:high_prob_guarantee}.
    With probability at least~$1-\beta$, the following hold true:
    \begin{itemize}
        \item {\bf Expectation}: The average-case performance $\Expect \big( \perfmetr[K](Z) \big)$ is upper-bounded by the optimal value of~\eqref{prob:dro-pep} with~$\alpha=1$, $J=1$, and $(\Aobj^1, \bobj^1, \cobj^1) = (\Aobj, \bobj, 0)$.
        \item $\cvar$: The~$\cvar$-performance $\cvar_\alpha (\perfmetr[K](Z))$ with~$\alpha\in(0, 1]$ is upper-bounded by the optimal value of~\eqref{prob:dro-pep} where~$J=2$,
        with~$(\Aobj^1, \bobj^1, \cobj^1) = \alpha^{-1} ( \Aobj, \bobj, \alpha - 1)$ and $(\Aobj^2, \bobj^2, \cobj^2) = (0, 0, 1)$.
    \end{itemize}
\end{corollary}

\begin{proof}
    From the finite-sample guarantee~\eqref{eq:high_prob_guarantee},
    the true distribution~$\prob_Z$ is an element of the Wasserstein ambiguity set~$\ambiset[\varepsilon]$ with probability at least~$1-\beta$.
    Therefore, $\Expect_{Z \sim \prob_Z} (\perfmetr[K](Z))$ and~$\cvar_\alpha (\perfmetr[K](Z))$ over~$\prob_Z$ are upper-bounded by the corresponding worst-case expectations over~$\probQ_Z \in \ambiset[\varepsilon]$, with probability at least~$1-\beta$.
    The proof then follows from applying Theorem~\ref{thm:conic-dro-pep} with the given choices of parameters and loss functions.
\end{proof}

\subsection{Proof of Theorem~\ref{thm:conic-dro-pep}}
\label{subsec:proof-conic-dro-pep}
We now derive the convex formulation stated in Theorem~\ref{thm:conic-dro-pep},
starting from
\begin{equation}
    \sup_{\probQ_Z \in \ambiset[\varepsilon]}\,
        \inf_{\ell \in \lclass}\,
        \Expect_{Z \sim \probQ_Z}\, \left( \ell(Z) \right)
    \;=\;
    \begin{array}[t]{ll}
        \sup
        &   \inf_{\ell \in \lclass}\, \Expect_{(G, F) \sim \probQ_Z}\, \left( \ell(G, F) \right) \\
        \text{s.t.}
        &   \supp\, \probQ_Z \subseteq \fclass \cap \xclass,
        \quad   \probQ_Z \in \ambiset[\varepsilon],
    \end{array}
    \label{prob:dro-bound}
\end{equation}
which is the right-hand side of~\eqref{eq:worst-case-exp}.
By adopting the proof strategies from~\citep[Theorem~4.2]{mohajerinesfahaniDatadrivenDistributionallyRobust2018} and~\citep[Appendix~A]{wangMeanRobustOptimization2024},
we first show that this problem is equivalent to the minimization problem~\eqref{prob:dro-pep}.
\begin{lemma}
    \label{lem:worst-case-exp-dual}
    Consider the setup of Theorem~\ref{thm:conic-dro-pep}.
    For every~$\varepsilon>0$, the worst-case expectation problem~\eqref{prob:dro-bound} and the minimization problem~\eqref{prob:dro-pep} have the same optimal value.
\end{lemma}

\begin{proof}
    Following the proof structure of~\citet[Theorem~4.2]{mohajerinesfahaniDatadrivenDistributionallyRobust2018}, we obtain the finite-dimensional dual problem to the worst-case expectation in~\eqref{prob:dro-bound} as
    \begin{align*}
        \min_{\lambda \in \reals_+,\, t \in \reals}
        \biggl\{
            \frac{1}{N} \sum_{i=1}^N
            \max_{Z_i \in \fclass \cap \xclass}
            \left(
                \ell(Z_i) - \lambda \norm{Z_i - \hZ_i}
            \right)
            + \lambda \varepsilon
        \biggr\},
    \end{align*}
    by interchanging~$\inf_{t \in \reals}$ and $\sup_{\probQ_Z \in \ambiset[\varepsilon]}$.
    Sion's minimax theorem justifies this interchange because the objective is convex and lower semicontinuous in~$t$, affine and continuous in~$\probQ_Z$, and the ambiguity set~$\ambiset[\varepsilon]$ is convex and compact on the compact support~$\fclass\cap\xclass$.
    Since~$\ell = \max_{1 \le j \le J} \ell^j$, the inner maximization splits into~$J$ maxima, one per branch, of~$\ell^j(Z_i) - \lambda \|Z_i - \hZ_i\|$ over~$Z_i \in \fclass \cap \xclass$.
    Applying Fenchel--Rockafellar duality~\citep[Theorem~15.23]{bauschkeConvexAnalysisMonotone2017} to the inner maximization for each~$j$ and introducing epigraph variables~$s_i$, we get
    \begin{equation*}
        \begin{array}[t]{ll}
            \textnormal{minimize}
            &   (1/N) \sum_{i=1}^N s_i \\
            \textnormal{subject to}
            &   [ -\ell^j ]^* ( \mu_i^j - \nu_i^j )
                + \sigma_{\fclass \cap \xclass} (\nu_i^j) 
                - \inner{\mu_i^j}{(\hG_i, \hF_i)}
                + \lambda \varepsilon \le s_i \\
            &   \norm[*]{\mu_i^j} \le \lambda,
                \hfill
                i=1,\dots,N,\, j=1,\dots,J,
        \end{array}
    \end{equation*}
    with variables~$s_i \in \reals$, $\lambda \in \reals_+$, $t\in\reals$, and~$\mu_i^j, \nu_i^j \in \symm^{K+2} \times \reals^{K+1}$.
    The Fenchel constraint qualification holds because~$-\ell^j$ and the mapping~$Z_i \mapsto \lambda \|Z_i - \hZ_i\|$ are finite and continuous on the ambient space, while~$\fclass \cap \xclass$ is nonempty, so for every fixed~$(\lambda, t)$ the dual representation for each branch is exact with attained minimum.
    Note that~$[-\ell^j]^*$ is a convex conjugate with respect to variables in~$\symm^{K+2} \times \reals^{K+1}$,~\ie,
    \begin{align*}
        [- \ell^j]^*(\mu_i^j - \nu_i^j)
        &=
            \begin{cases}
                \cobj^j t & \textnormal{if } \mu_i^j - \nu_i^j + (\Aobj^j, \bobj^j) = 0, \\
                \infty & \textnormal{otherwise}.
            \end{cases}
    \end{align*}
    The support function~$\sigma_{\fclass \cap \xclass}$ of~$\fclass \cap \xclass$ as in Assumption~\ref{assumption:interpolation}
    can be written as
    \begin{align*}
        \sigma_{\fclass \cap \xclass} (\nu_i^j)
&=
            \begin{array}[t]{llll}
                \max & \inner{\nu_i^j}{(G, F)}
                &   \textnormal{s.t.}
                &   \linear(G, F) \in \reals_+^{|M|},
                \\ & & &
                    - \inner{(A_0, b_0)}{(G, F)} - c_0 \in \reals_+
            \end{array} \\
        &=
            \begin{array}[t]{llll}
                \min & - c_0 \tau_i^j
                &   \textnormal{s.t.}
                &   - \linear^*(y_i^j) + \tau_i^j (A_0, b_0) - \nu_i^j \in \symm_+^{K+2} \times \{0\},
            \end{array}
    \end{align*}
    with variables~$(G, F) \in \symm_+^{K+2} \times \reals^{K+1}$, $\tau_i^j \in \reals_+$, and~$y_i^j \in \reals_+^{|M|}$.
    For each fixed~$\nu_i^j$, the support-function maximization is an instance of the performance estimation problem with an objective linear in~$(G, F)$.
    The function class and nondegeneracy condition $\eta_k^k \neq 0$ of Theorem~\ref{thm:conic-dro-pep} allow us to apply~\citep[Theorem~6]{taylorSmoothStronglyConvex2017},
    so its conic dual attains the same optimal value, which gives the last equality.
    The conjugate constraint gives~$\nu_i^j = \mu_i^j + (\Aobj^j, \bobj^j)$, so setting~$\mu_i^j = (X_i^j, Y_i^j)$ eliminates~$\nu_i^j$ from the conic constraint and gives
    \begin{equation*}
        - \linear^*(y_i^j) - (X_i^j, Y_i^j)
        + \tau_i^j (A_0, b_0) - (\Aobj^j, \bobj^j)
        \in \symm_+^{K+2} \times \{0\}.
    \end{equation*}
    Substituting this constraint and the conjugate value~$\cobj^j t$ into the epigraph problem for each branch gives exactly~\eqref{prob:dro-pep}.
\end{proof}

Next, we dualize the minimization problem to obtain the maximization problem.
\begin{lemma}
    \label{lem:worst-case-exp-double-dual}
    Consider the setup of Theorem~\ref{thm:conic-dro-pep}.
    For every~$\varepsilon>0$, strong duality holds between the minimization problem~\eqref{prob:dro-pep} and the maximization problem of Theorem~\ref{thm:conic-dro-pep}, whose dual, or maximization, variables are~$(G_i^j, F_i^j) \in \fclass \cap \xclass$ and~$(\pi_i^j) \in \Gamma$.
\end{lemma}

\begin{proof}
    First, we write the dual-norm constraint~$\|(X_i^j, Y_i^j)\|_* \le \lambda$ of~\eqref{prob:dro-pep} as the second-order cone constraint~$\big( \lambda_i^j, (X_i^j, Y_i^j) \big) \in \soc$ with auxiliary variables $\lambda_i^j = \lambda$ for~$i=1,\dots,N$ and~$j=1,\dots,J$.
    Then the Lagrangian of this problem is
    \begin{align*}
        &\lagrangian \left(
            s_i, \lambda, t, y_i^j, \tau_i^j, \big( \lambda_i^j, (X_i^j, Y_i^j) \big),\;
            \pi_i^j, (G_i^j, F_i^j), \varepsilon_i^j
        \right) \\
&=
            {\textstyle \sum_{i=1}^N} s_i \big( 1/N - {\textstyle \sum_{j=1}^J} \pi_i^j \big)
            + \lambda {\textstyle \sum_{i,j}} ( \pi_i^j \varepsilon - \varepsilon_i^j )
            + t\, {\textstyle \sum_{j=1}^J} ( {\textstyle \sum_{i}} \pi_i^j ) \cobj^j
        \\&\quad
            + {\textstyle \sum_{i,j}} \Big[ (y_i^j)^T \linear(G_i^j, F_i^j)
            - \tau_i^j \big( \inner{(A_0, b_0)}{(G_i^j, F_i^j)} + \pi_i^j c_0 \big)
            + \lambda_i^j \varepsilon_i^j
        \\&\qquad\qquad\quad
            + \inner{(X_i^j, Y_i^j)}{(G_i^j, F_i^j) - \pi_i^j (\hG_i, \hF_i)}
            + \inner{(\Aobj^j, \bobj^j)}{(G_i^j, F_i^j)} \Big],
    \end{align*}
    with dual variables~$\pi_i^j \in \reals_+$, $(G_i^j, F_i^j) \in \symm_+^{K+2} \times \reals^{K+1}$, and $\varepsilon_i^j \in \reals$ for~$i=1,\dots,N$ and~$j=1,\dots,J$.
    We minimize this over the primal variables to obtain the dual problem; 
    in particular, minimizing over~$\big( \lambda_i^j, (X_i^j, Y_i^j) \big) \in \soc$ forces~$\big( \varepsilon_i^j,\, (G_i^j, F_i^j) - \pi_i^j (\hG_i, \hF_i) \big) \in \soc$, and~$\varepsilon_i^j$ is then eliminated using~$\big\| (G_i^j, F_i^j) - \pi_i^j (\hG_i, \hF_i) \big\| \le \varepsilon_i^j$.
    \begin{equation*}
        \begin{array}[t]{llr}
            \textnormal{maximize}
            &   \sum_{i=1}^N \sum_{j=1}^J \big\langle (\Aobj^j, \bobj^j),\, (G_i^j, F_i^j) \big\rangle \\
            \textnormal{subject to}
                &   \sum_{j=1}^J \pi_i^j = 1/N,  &   i=1,\dots,N \\
                &   \sum_{j=1}^J \left( \sum_{i=1}^N \pi_i^j \right) \cobj^j = 0 \\
                &   \sum_{i=1}^N \sum_{j=1}^J \big\| (G_i^j, F_i^j) - \pi_i^j (\hG_i, \hF_i) \big\| \le \varepsilon \\
                &   \linear(G_i^j, F_i^j) \in \reals_+^{|M|},  & i=1,\dots,N,\; j=1,\dots,J, \\
                &   \big\langle (A_0, b_0),\, (G_i^j, F_i^j) \big\rangle + \pi_i^j c_0 \le 0,  & i=1,\dots,N,\; j=1,\dots,J,
        \end{array}
    \end{equation*}
    with variables~$(G_i^j, F_i^j) \in \symm_+^{K+2} \times \reals^{K+1}$ and $\pi_i^j \in \reals_+$.
    Applying the change of variables~$(G_i^j, F_i^j) \leftarrow \pi_i^j (G_i^j, F_i^j)$ scales the objective and the norm constraint by~$\pi_i^j$, while the remaining constraints become~$\pi_i^j \linear(G_i^j, F_i^j) \in \reals_+^{|M|}$ and~$\pi_i^j \big( \inner{(A_0, b_0)}{(G_i^j, F_i^j)} + c_0 \big) \le 0$.
    For~$\pi_i^j > 0$ these are equivalent to~$(G_i^j, F_i^j) \in \fclass \cap \xclass$, while for~$\pi_i^j = 0$ the variable~$(G_i^j, F_i^j)$ does not affect the objective and~$(\hG_i, \hF_i) \in \fclass \cap \xclass$ is feasible.
    We may therefore impose~$(G_i^j, F_i^j) \in \fclass \cap \xclass$ throughout, obtaining the maximization problem of Theorem~\ref{thm:conic-dro-pep} with variables~$(G_i^j, F_i^j) \in \fclass \cap \xclass$ and~$(\pi_i^j) \in \Gamma$.

    Let~$y_\star^j \in \reals_+^{|M|}$ and~$\tau_\star^j \in \reals_+$ be a feasible solution of the original dual PEP problem, so~$-\linear^*(y_\star^j)+\tau_\star^j(A_0,b_0)-(\Aobj^j,\bobj^j)=(S_\star^j,0)$ for some~$S_\star^j \succeq 0$.
    Fix~$\delta,\omega>0$, define~$(D_G,D_F)=-\linear^*(\ones)+(A_0,b_0)$, and set~$y_i^j=y_\star^j+\delta\ones$, $\tau_i^j=\tau_\star^j+\delta$, and~$(X_i^j,Y_i^j)=\delta(D_G,D_F)-(\omega I,0)$.
    Then~$y_i^j>0$, $\tau_i^j>0$, and
    \begin{equation*}
        -\linear^*(y_i^j)-(X_i^j,Y_i^j)
        +\tau_i^j(A_0,b_0)-(\Aobj^j,\bobj^j)
        =(S_\star^j+\omega I,0)
        \in\symm_{++}^{K+2}\times\{0\}.
    \end{equation*}
    Taking~$\lambda>\max_{i,j}\|(X_i^j,Y_i^j)\|_*$,~$t=0$, and~$s_i$ sufficiently large yields a Slater point for~\eqref{prob:dro-pep} and hence strong duality.
\end{proof}
Lemmas~\ref{lem:worst-case-exp-dual} and~\ref{lem:worst-case-exp-double-dual}
establish the two equivalences in Theorem~\ref{thm:conic-dro-pep},
the latter with strong duality, and thus complete its proof.

\subsection{Bridging worst-case and average-case convergence rates}
\label{subsec:rates_interpolation}

The radius~$\varepsilon$ of the Wasserstein ambiguity set~$\ambiset[\varepsilon]$ controls the trade-off between \emph{a priori} and \emph{a posteriori} information.
As~$\varepsilon$ grows, the upper bounds in Corollary~\ref{cor:conic-dro-pep} approach the worst-case convergence bound, which only uses \emph{a priori} information, \ie, function classes and initial conditions.
As~$\varepsilon$ shrinks, the bounds from Corollary~\ref{cor:conic-dro-pep} approach the \emph{a posteriori} empirical average and~$\cvar$ losses.
The following theorem characterizes the behavior of the problem~\eqref{prob:dro-pep} in these extreme cases.

\begin{theorem}
    \label{thm:pep-interpolation}
    Consider the setup of Corollary~\ref{cor:conic-dro-pep}.
    \begin{itemize}
        \item
        When~$\varepsilon \to 0^+$,
        the optimal values of~\eqref{prob:dro-pep} converge respectively to the in-sample average-case performance~$(1/N) \sum_{i=1}^N \perfmetr[K] (\hZ_i)$ and the in-sample $\cvar$-performance~$\cvar_\alpha \big( \perfmetr[K](\hZ) \big)$.

        \item
        There exists $\bar{\varepsilon} > 0$ such that for any~$\varepsilon \ge \bar{\varepsilon}$, the optimal values of~\eqref{prob:dro-pep}
        all equal the worst-case performance.
    \end{itemize}
\end{theorem}

\begin{proof}
    First, the case of~$\varepsilon \to 0^+$ follows from~\citep[Proposition~8.5]{kuhnDistributionallyRobustOptimization2025}; the optimal value of~\eqref{prob:dro-pep} converges to the in-sample loss value.
    Now consider the case of large~$\varepsilon$ values.
    By Lemma~\ref{lem:domain_size}, $\fclass \cap \xclass$ is compact.
    Then there exists a worst-case instance~$\bar{Z} \in \fclass \cap \xclass$ such that~$\perfmetr[K](\bar{Z}) = \sup_{Z \in \fclass \cap \xclass} \perfmetr[K](Z)$.
    Let~$\bar{\varepsilon} = \max_{1 \le i \le N} \|\bar{Z} - \hZ_i\|$.
    Then for any~$\varepsilon \ge \bar{\varepsilon}$, setting~$Z_i^j = \bar{Z}$ for every~$i=1,\dots,N$ and~$j=1,\dots,J$ is feasible in the maximization problem of Theorem~\ref{thm:conic-dro-pep}, since
    \begin{equation*}
        {\textstyle \sum_{i=1}^N \sum_{j=1}^J} \pi_i^j \|\bar{Z} - \hZ_i\|
        = (1/N) {\textstyle \sum_{i=1}^N} \| \bar{Z} - \hZ_i \|
        \le \bar{\varepsilon}
        \le \varepsilon,
    \end{equation*}
    for any~$(\pi_i^j) \in \Gamma$.
    Therefore, $\perfmetr[K](\bar{Z})$ simultaneously lower-bounds and upper-bounds the worst-case expectations in Corollary~\ref{cor:conic-dro-pep}, which concludes the proof.
\end{proof}

\section{Numerical experiments}\label{sec:experiments}
We now analyze the probabilistic performance of various first-order methods in unconstrained quadratic optimization, logistic regression, and Lasso.
All code to reproduce our experiments is available at
\begin{center}
    \url{https://github.com/stellatogrp/dro_pep}.
\end{center}
We solve the minimization problems in Corollary~\ref{cor:conic-dro-pep}.
We implement all examples in Python 3.13, formulate~\eqref{prob:dro-pep} with JAX 0.9.0, and solve the resulting problems with Clarabel 0.11 under default settings~\citep{goulartClarabelInteriorpointSolver2026}.
We compute worst-case convergence bounds with PEPit 0.5.1~\citep{goujaudPEPitComputerassistedWorstcase2024}.

\paragraph{Parameter selection}
For each experiment, we generate $N=100$ i.i.d.\ samples, each consisting of function realizations and corresponding iterate trajectories, to form the~\eqref{prob:dro-pep} problem.
Section~\ref{subsec:DRO} shows how the radius~$\varepsilon$ controls the finite-sample guarantee when the support constants are valid for the full instance distribution.
In the experiments, we instead calibrate~$\varepsilon$ empirically to target 95\% validation coverage.
We draw 100 independent repetitions of the sampling experiment, each consisting of a fresh batch of problem instances (100 for the quadratic experiment, 200 for logistic regression and Lasso), and compute the empirical statistic (mean or~$\cvar$) of each repetition.
For each~$K$, we select the smallest tested~$\varepsilon$ whose certificate upper-bounds the 95\textsuperscript{th} percentile of these statistics.
This is empirical calibration on validation batches independent of the~$N$ training instances.
It does not evaluate coverage on an additional held-out set, and it is distinct from using a distribution-free radius~$\varepsilon_N(\beta)$ in~\eqref{eq:high_prob_guarantee}.
Likewise, support constants estimated from reference instances below are empirical calibrations, not uniform bounds on distributions with unbounded Gaussian components.

\paragraph{Numerical results}
Each experiment compares the following across iterations:
\begin{itemize}
    \item \textbf{Worst-case convergence bound from PEP.}
    We use the PEPit toolbox to solve the corresponding PEP for worst-case convergence bounds.

    \item \textbf{DRO objective values.}
    We solve the~\eqref{prob:dro-pep} problems in Corollary~\ref{cor:conic-dro-pep} to obtain distributionally robust expectation and~$\cvar$ values for the chosen performance metric, using the empirically calibrated radii.
    In all plots, the solid lines show the optimal objective value of~\eqref{prob:dro-pep}.
\end{itemize}

\subsection{Unconstrained quadratic minimization}
\label{subsec:quad_experiment}
For this experiment, we restrict ourselves to~$L$-smooth convex \emph{quadratic} functions
\begin{equation*}
    f(x) = (1/2) (x - x^\star)^T Q(x - x^\star),
\end{equation*}
where~$Q \in \symm_+^d$ is a positive semidefinite matrix whose eigenvalues are contained in~$[0, L]$ and $x^\star$ is an optimal solution.
For scalability reasons, we solve~\eqref{prob:dro-pep} with the function class~$\mathcal{F}_{0,L}$ of~$L$-smooth convex functions,
while computing the worst-case performance in PEPit with the tightest function class~$\mathcal{Q}_{0, L}$ of convex quadratics~\citep{bousselmiInterpolationConditionsLinear2024}.

We compare gradient descent~\eqref{eq:GD} and Nesterov's fast gradient method~\eqref{eq:FGM}, both with step size~$\eta = 1/L$.
We run FGM with the momentum coefficient~$\gamma_k = (\theta_k - 1) / \theta_{k+1}$, where~$\{\theta_k\}$ is defined by~$\theta_{k+1} = (1 + \sqrt{1 + 4 \theta_k^2})/2$ from~$\theta_0 = 1$~\citep{nesterovMethodUnconstrainedConvex1983}.
This differs from the canonical momentum~$k/(k+3)$ in Proposition~\ref{prop:ls_rate}; its FGM rates are therefore reference curves, not theorem predictions for the implemented recurrence.
The sampled trajectories, the worst-case PEP bounds, and the DRO certificates all use these coefficients.

\paragraph{Problem distribution}
We sample values of~$Q$, whose eigenvalues are all in~$[0, L] \subseteq \reals$, to match existing results on average-case analysis~\citep{cunhaOnlyTailsMatter2022,paquetteHaltingTimePredictable2023,scieurUniversalAverageCaseOptimality2020} and our rates in Proposition~\ref{prop:ls_rate}.
Define~$A \in \reals^{d \times d}$ whose entries are i.i.d.\ Gaussian~$\mathcal{N}(0, L/4)$.
In the regime~$d \to \infty$, the empirical spectral distribution~$d^{-1} \sum_{i=1}^d \delta_{\lambda_i}$ of~$Q = (1/d)A^T A$ converges to the \emph{Marchenko--Pastur (MP) distribution}~\citep{marcenkoDistributionEigenvaluesSets1967}.
Its density is~$\mathrm d \mu_{\mathrm{MP}} (\lambda) = 2 \sqrt{(L-\lambda)\lambda} / (\pi L \lambda) \mathbf{1}_{[0,L]}(\lambda) \, \mathrm d\lambda$, where $\mathbf{1}_{[0, L]}$ is the indicator of~$[0, L]$.
Note that this distribution corresponds to the case of~$a = -1/2$ in Proposition~\ref{prop:ls_rate}.
We use sufficiently large~$d$ and rejection sampling to ensure all sample values of~$Q$ have eigenvalues in our predefined range~$[0, L]$.
We sample the entries of~$x^\star$ i.i.d.\ from~$\mathcal{N}(0,1)$ as well.

\begin{figure}[t]
    \centering
        \includegraphics[width=\textwidth]{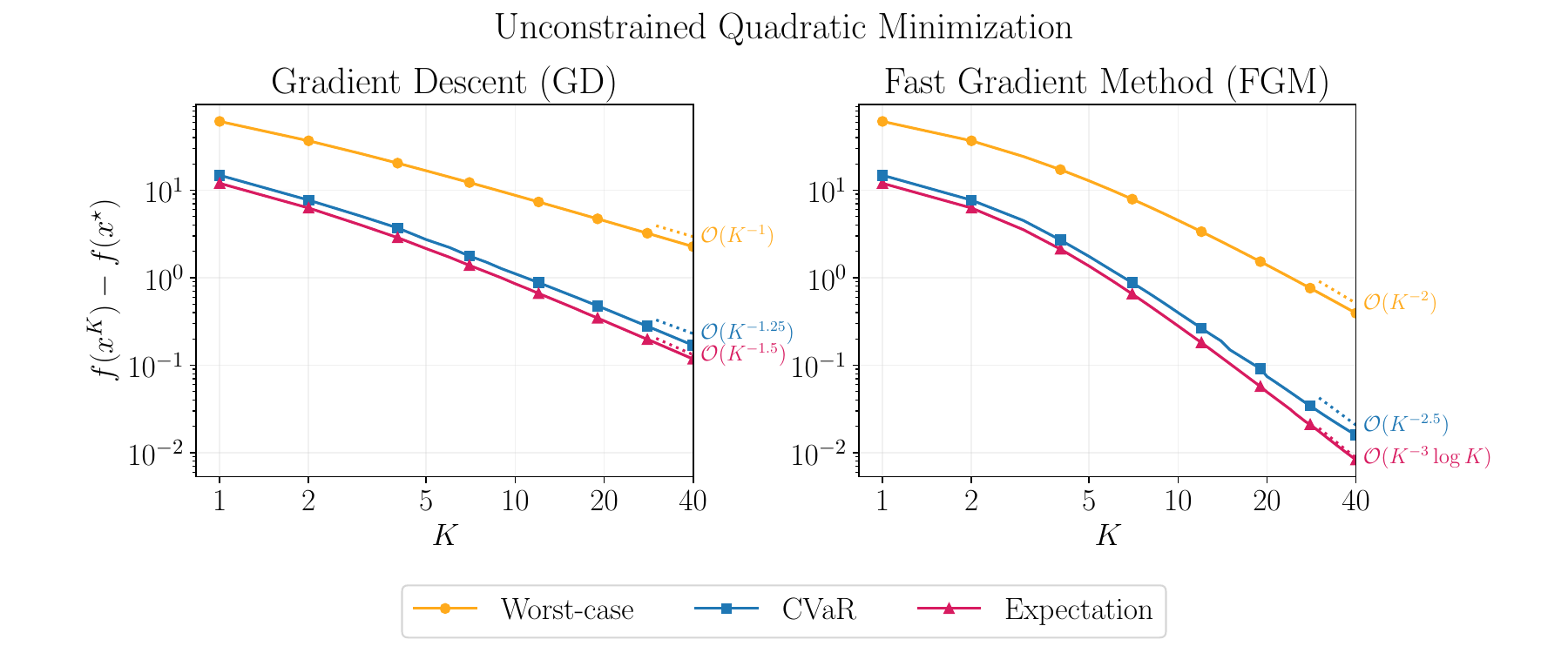}
    \caption{
        Function-value gaps for GD and FGM on the unconstrained quadratic instances.
        The curves compare the worst-case PEP bound with empirically calibrated DRO certificates for the expectation and~$\cvar_{0.05}$; the labeled rate lines are theoretical references.
    }
    \label{fig:quad_mp_comp_bound}
\end{figure}

\paragraph{Problem setup}
All runs start from the origin~$x^0 = 0$.
We set~$d=300$,~$L=25$, and~$\eta=1/L$, and calibrate the initial-condition radius in~$\|x^0 - x^\star\| \le r$ as the largest distance observed over~$200$ reference instances, giving~$r = 19.173$.
We use the function-value gap~$f(x^K) - f^\star$ as the performance metric for~$K = 1, \dots, 40$, and compare our bounds with the GD rates and the canonical-momentum FGM reference rates of Proposition~\ref{prop:ls_rate}.
For the~$\cvar$ experiment, we set $\alpha=0.05$.
The calibration grid consists of 23 values of~$\varepsilon$ spanning~$[10^{-5},10]$.
It combines a logarithmically spaced grid on~$[10^{-1},10]$ with smaller values that ensure the selected~$\varepsilon$ is not the smallest candidate.

\paragraph{Comparison with known probabilistic rates}
As discussed in Section~\ref{subsec:risk_measure}, probabilistic analysis can be challenging as it requires full knowledge of the underlying distribution~$\prob$.
The analysis of~\citet{cunhaOnlyTailsMatter2022,paquetteHaltingTimePredictable2023} derives the exact average-case rates~$\Theta(K^{-3/2})$ for GD and $\Theta(K^{-3} \log K)$ for canonical-momentum FGM in terms of the suboptimality~$f(x^K) - f^\star$.
Our analysis of Proposition~\ref{prop:ls_rate} ($a = -1/2$) further specifies the~$\cvar$-performance upper bound:
the bound is~$O(K^{-1.25})$ for GD and~$O(K^{-2.5})$ for canonical-momentum FGM.

Figure~\ref{fig:quad_mp_comp_bound} shows a GD certificate consistent with its theoretical rate, without fitting the spectral distribution explicitly; the canonical-momentum FGM curve is only a qualitative reference for the accelerated run.

\subsection{Logistic regression}
\label{subsec:logreg_experiment}
Consider a binary classification problem with data~$\{(a_i, b_i)\}_{i=1}^m$, where~$a_i \in \reals^d$ are feature vectors and~$b_i \in \{0, 1\}$ are labels.
The goal is to learn a vector~$x \in \reals^d$ that linearly models the log-odds of a feature vector belonging to the positive class.
The associated objective is
\begin{equation}\label{eq:logreg}
    f(x) = -\frac{1}{m}\sum_{i=1}^m \Big( b_i \log \sigma(a_i^T x) + (1-b_i)\log\big(1 - \sigma(a_i^T x)\big) \Big),
\end{equation}
where~$\sigma(z) = 1/(1 + e^{-z})$ is the sigmoid function~\citep[Chapter 7]{boydConvexOptimization2004}.
Let~$A \in \reals^{m \times d}$ be the matrix of feature vectors stacked as rows.
The function~$f$ is convex and~$L$-smooth with~$L = \lambdamax(A^T A)/(4m)$~\citep[Theorem 5.12]{beckFirstOrderMethodsOptimization2017}.
We do not add a regularizer, so the strong convexity parameter is~$\mu = 0$.

We again compare gradient descent~\eqref{eq:GD} and Nesterov's fast gradient method~\eqref{eq:FGM} with the same momentum coefficients as in Section~\ref{subsec:quad_experiment}.
Since GD converges on smooth convex functions for any step size~$\eta \in (0, 2/L)$, we run it with the larger step size~$\eta = 1.9/L$, selected by tuning on sampled instances, while FGM uses~$\eta = 1/L$.
The worst-case PEP baselines are computed at the same step sizes.

\paragraph{Problem distribution}
We construct a distribution over logistic regression instances from the \texttt{german.numer} credit-scoring dataset in the LIBSVM repository~\citep{changLIBSVMLibrarySupport2011}, which contains~$1000$ examples with~$24$ numerical features.
We standardize each feature to zero mean and unit variance over the full dataset and append a constant intercept feature, so that~$d = 25$.
Each problem instance is the objective~\eqref{eq:logreg} on an independent uniformly random subsample of~$m = 300$ examples, so the instances are i.i.d.\ samples from this subsampling distribution.
This models the common setting where the same method repeatedly solves problems on minibatches or resampled versions of a dataset.
We verified numerically that every sampled and reference instance had a finite minimizer, as required for the unregularized model.

\paragraph{Problem setup}
All runs start at the origin~$x^0 = 0$.
As in Section~\ref{subsec:quad_experiment}, we use the largest constants observed over~$200$ reference instances for the PEP comparison: the smoothness~$L = \max_i \lambdamax(A_i^T A_i)/(4m) = 0.770$ and the initial distance radius~$r = \max_i \|x^0 - x^\star_i\|_2 = 8.14$.
We use the squared gradient norm~$\|\nabla f(x^K)\|^2$ as the performance metric for~$K = 1, \dots, 30$.
The calibration grid consists of 16 logarithmically spaced values~$\varepsilon \in [10^{-8}, 10^{-1/2}]$, reaching radii small enough that the selected~$\varepsilon$ is never the smallest candidate.
Moreover, $\alpha = 0.01$ in the~$\cvar$ experiment.

\begin{figure}[t]
    \centering
        \includegraphics[width=0.95\textwidth]{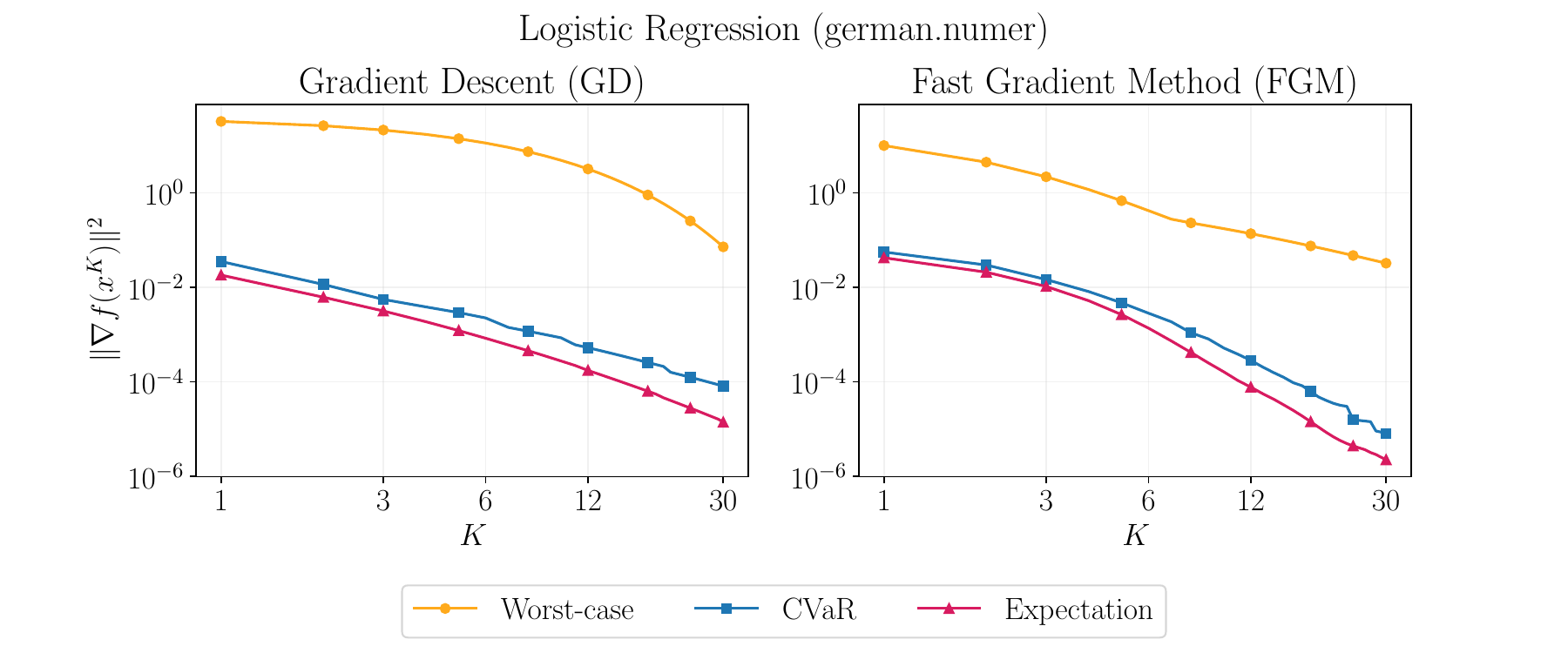}
    \caption{
        Squared-gradient performance for GD and FGM on logistic regression problems from subsamples of the \texttt{german.numer} dataset~\citep{changLIBSVMLibrarySupport2011}.
        The curves compare the worst-case PEP bound with empirically calibrated DRO certificates for the expectation and~$\cvar_{0.01}$.
    }
    \label{fig:logreg_obj_val_bound_plots}
\end{figure}

\paragraph{Results}
Figure~\ref{fig:logreg_obj_val_bound_plots} shows the convergence bounds for~$\perfmetr[K](f, x^0) = \|\nabla f(x^K)\|^2$ as a function of~$K$.
The worst-case PEP bound sits two to four orders of magnitude above the data-driven certificates at every~$K$.
At~$x^0 = 0$, the vector~$\sigma(Ax^0)-b$ has squared norm~$m/4$.
Thus, every logistic instance satisfies~$\|\nabla f(x^0)\|^2 \le \lambdamax(A^T A)/(4m) \le L = 0.770$.
Functions in~$\mathcal{F}_{0,L}$ with the same constants~$(L,r)$ can instead attain initial squared gradient norms up to~$L^2r^2 \approx 39$.
The PEP baseline must cover these more adversarial functions, whereas the data-driven certificates adapt to the smaller gradients of the logistic instances.
The GD worst-case curve decays geometrically over the plotted range, consistent with the two-branch structure of exact gradient-norm bounds for constant-step GD~\citep{rotaruExactWorstcaseConvergence2024}.
For the distance-bounded PEP used here, the step size~$\eta = 1.9/L$ keeps the geometric branch~$|1 - \eta L|^{2K}$ active until~$K \approx 40$.
The expectation certificates closely track the empirical means, while the gap between the~$\cvar$ and expectation certificates widens as~$K$ grows.
Average-case analysis therefore increasingly understates the residual performance on the hardest instances, making risk-aware certificates more informative at later iterations.

\subsection{Lasso}\label{subsec:lasso_experiment}
Consider an~$\ell_1$-regularized least-squares problem, often referred to as the Lasso~\citep{tibshiraniRegressionShrinkageSelection1996}.
The objective of this problem is
\begin{equation}\label{eq:lasso}
    f(x) = (1/2) \norm[2]{Ax - b}^2 + \lambda \norm[1]{x},
\end{equation}
with~$A \in \reals^{m \times d}$, $b \in \reals^m$, and regularization parameter~$\lambda > 0$.

As the objective function~\eqref{eq:lasso} is nondifferentiable, we analyze two proximal methods~\citep{parikhProximalAlgorithms2014}: the iterative soft-thresholding algorithm (ISTA) and the fast iterative soft-thresholding algorithm (FISTA)~\citep{beckFastIterativeShrinkageThresholding2009}.
Given parameter~$\delta > 0$, the soft-thresholding function is given by elementwise application of~$\mathcal{T}_\delta(v) = \min \big\{ v + \delta,\, \max \{ v - \delta,\, 0 \} \big\}$.
The ISTA with step size~$\eta > 0$ is defined by
\begin{equation}\label{eq:ista}\tag{ISTA}
    x^{k+1} = \mathcal{T}_{\lambda \eta}(x^k - \eta A^T(Ax^{k} - b)),
    \qquad k=0,1,\dots.
\end{equation}
The FISTA iterations with step size~$\eta > 0$ are given by
\begin{equation}\label{eq:fista}\tag{FISTA}
    \begin{array}{ll}
        x^{k+1} & = \mathcal{T}_{\lambda \eta}(y^k - \eta A^T(Ay^{k} - b))\\
        y^{k+1} & = x^{k+1} + \gamma_k(x^{k+1} - x^k),
    \end{array}
    \qquad k=0,1,\dots,
\end{equation}
with initialization~$x^0 = y^0$ and the momentum coefficient~$\gamma_k = (\theta_k - 1) / \theta_{k+1}$ used in the FGM experiments of Section~\ref{subsec:quad_experiment}.

\paragraph{Problem distribution}
We use a sparse coding example that recovers a sparse vector~$\tilde{x} \in \reals^d$ from noisy linear measurements~$b = A\tilde{x} + \xi \in \reals^m$~\citep{chenTheoreticalLinearConvergence2018}.
Here~$A \in \reals^{m \times d}$ is a known dictionary matrix, and the noise satisfies~$\xi_i \sim \mathcal{N}(0, \sigma_\xi^2)$.
We sample one value of~$A$ for the entire experiment, both training and testing, and construct a distribution of instances by sampling values of~$\tilde{x}$ and~$\xi$ to form~$b$.
We generate~$A$ with a controlled spectrum: we draw orthogonal matrices~$U \in \reals^{m \times m}$ and~$V \in \reals^{d \times d}$, and set~$A = U \Sigma V^T$, where~$\Sigma \in \reals^{m \times d}$ is rectangular diagonal with singular values~$s$ chosen so that~$A^T A$ has a two-level eigenvalue spectrum on its range: a cluster of value~$q L$ together with a single top spike at~$L$, normalized so that the largest eigenvalue is exactly~$L=1$.
After generating~$A$, we use the value of the smoothness parameter~$L$ (the top eigenvalue of~$A^T A$) and an experimentally computed initial distance radius~$r$ for the PEP comparison; since~$m < d$, the strong convexity parameter is~$\mu = 0$.

\paragraph{Problem setup}
We set~$m = 100$,~$d = 150$,~$\lambda=3\times10^{-4}$,~$q = 0.08$,~$\sigma_\xi = 10^{-2}$, and~$\eta = 1/L$.
Each instance is formed by sampling the coordinates independently, with~$\tilde{x}_i \sim \mathcal{N}(0,1)$ with probability~$0.2$ and~$\tilde{x}_i = 0$ otherwise, then setting~$b = A\tilde{x} + \xi$.
All runs start from the origin~$x^0 = 0$, so the initial condition~$\|x^0 - x^\star\| \le r$ is enforced against the per-instance solution~$x^\star(b)$; the largest distance over~$200$ reference instances gives~$r=7.482$.
We run both methods for~$K = 1, \dots, 25$.
For the~$\cvar$ experiment, we again set $\alpha=0.05$.
The calibration grid consists of 19 values of~$\varepsilon$ spanning~$[10^{-5}, 10^{-1}]$, refining a linearly spaced base grid on~$[10^{-3}, 10^{-1}]$ so that the selected~$\varepsilon$ always lies in the grid interior.

\begin{figure}[ht]
    \centering
        \includegraphics[width=0.95\textwidth]{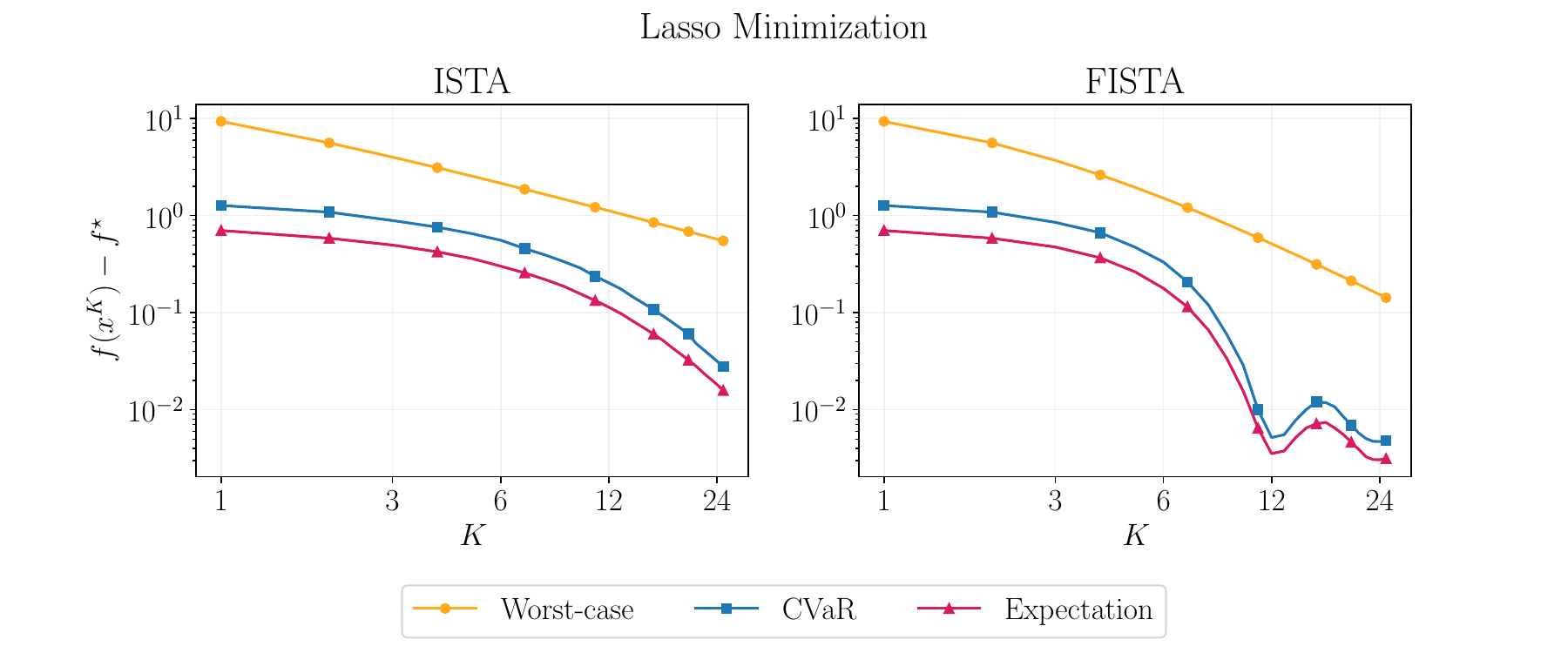}
    \caption{
        Function-value gaps for ISTA and FISTA on the Lasso instances.
        The curves compare the worst-case PEP bound with empirically calibrated DRO certificates for the expectation and~$\cvar_{0.05}$.
    }
    \label{fig:lasso_bound_plots}
\end{figure}

\paragraph{Results}
Figure~\ref{fig:lasso_bound_plots} shows the convergence bounds for~$\perfmetr[K](f, x^0) = f(x^K) - f^\star$ as a function of~$K$.
For this distribution, we observe that ISTA exhibits a generally sublinear convergence trend that is monotonically decreasing in~$K$.
In contrast, FISTA's behavior exhibits the rippling effect of the momentum term, and our DRO certificates are able to capture this phenomenon.

\section{Conclusion}
In this paper, we presented a probabilistic performance estimation framework for deterministic first-order methods applied to problems drawn from a distribution.
By combining the performance estimation framework with distributionally robust optimization, we derived data-driven probabilistic performance guarantees from iterate trajectories on finitely many sampled instances.
We highlighted the effectiveness of our method by analyzing the setup where both a priori and a posteriori information are necessary to obtain tight convergence rates.
In particular, our framework extends PEP from worst-case analysis to data-driven guarantees for the expectation and conditional value-at-risk.

A key limitation of our approach is the scalability of~\eqref{prob:dro-pep}, as its dimension is proportional to the number of samples and the square of the number of iterations.
One way to address this issue is to adopt scenario reduction techniques, such as mean robust optimization~\citep{wangMeanRobustOptimization2024}, that cluster data points to reduce the size of the DRO problem.
A further promising direction is algorithm design. Rather than certifying a fixed method, one can optimize the method parameters directly against the distributionally robust performance bound.
Preliminary results explore this direction by learning step sizes that minimize the DRO performance bounds~\citep{ranjanDistributionallyRobustLearning2026}.
Deriving optimal closed-form step sizes for structured distributions of instances is an interesting direction for future work.

\ifpreprint
\subsection*{Acknowledgements}
\myack
\fi

\ifpreprint\else
\acks{\myack}
\fi

\appendix

\section{Composite optimization}
\label{appendix:composite}
Consider the problem
\begin{equation*}
    \begin{array}[t]{ll}
        \text{minimize} & f(x) + h(x),
    \end{array}
\end{equation*}
where~$x \in \reals^d$ is the optimization variable, $f \colon \reals^d \to \reals$ is convex, and~$h \colon \reals^d \to \reals \cup \{\infty\}$ is proper, lower semicontinuous, and convex.

\begin{assumption}\label{asm:composite_h}
    We consider composite objectives~$f + h$ with an optimal solution~$x^\star \in \argmin_x (f + h)$ and finite optimum~$(f+h)^\star = f(x^\star) + h(x^\star) > -\infty$,
    with stationarity condition~$0 \in \nabla f(x^\star) + \partial h(x^\star)$ and~$s^\star = - \nabla f(x^\star) \in \partial h(x^\star)$.
    The smooth part~$f \in \fclass$ is continuously differentiable with~$L$-Lipschitz gradient, \eg,~$\fclass = \mathcal{F}_{\mu, L}$, $\mathcal{F}_{0, L}$, or~$\mathcal{Q}_{\mu, L}$.
The nonsmooth part~$h \in \fclassprox$ may be nondifferentiable, \eg, an indicator function~$\iota_C$ of a closed convex set~$C$ or the~$\ell_1$ norm.
    Without loss of generality, let~$x^\star = 0$ and~$(f+h)^\star = 0$.
\end{assumption}

\paragraph{Algorithm}
We say that~$\algo \colon (f, h, x^0) \mapsto \{x^k\}_{k=0,1,\dots}$ is a \emph{fixed-step proximal gradient method} if it extrapolates past iterates, applies a fixed-step smooth half-step, and then takes a proximal step~\citep{parikhProximalAlgorithms2014}:
\begin{equation*}
    x^{k+1} = \prox_{\gamma^{k+1} h}(x^{k+1/2}),
    \qquad
    x^{k+1/2} = y^k - {\textstyle \sum_{i=0}^{k}} \eta_k^i\, \nabla f(y^i),
    \qquad
    y^k = {\textstyle \sum_{j=0}^{k}} \beta_k^j\, x^j,
\end{equation*}
for~$k = 0, 1, \dots$, with proximal step size~$\gamma^{k+1}$, gradient update step sizes~$\{\eta_k^i\}_{i=0}^{k}$, and extrapolation coefficients~$\{\beta_k^j\}_{j=0}^{k}$, which we assume to satisfy the following.
\begin{assumption}\label{assumption:stepsize}
    The proximal step sizes, gradient update step sizes, and extrapolation coefficients are fixed independently of the problem instance~$(f,h,x^0)$.
    The proximal step sizes satisfy~$\gamma^{k+1}>0$, and the coefficients satisfy
    \begin{equation*}
        {\textstyle \sum_{i=0}^{k}} \eta_k^i = \gamma^{k+1},
        \qquad {\textstyle \sum_{j=0}^k} \beta_k^j = 1,
        \qquad k = 0, 1, \dots.
    \end{equation*}
    When~$h \equiv 0$, the proximal map is the identity and the proximal step size has no effect on the iterates.
    In this case, we drop the coupling condition~$\sum_{i=0}^{k} \eta_k^i = \gamma^{k+1}$ and allow arbitrary fixed gradient update coefficients~$\{\eta_k^i\}_{i=0}^k$.
\end{assumption}
Momentum-free methods set~$y^k=x^k$, as in~\eqref{eq:ista}, whereas FISTA~\eqref{eq:fista} extrapolates the two most recent iterates.
Let~$g^i=\nabla f(y^i)$ for~$i\ge0$ and~$s^i=(x^{i-1/2}-x^i)/\gamma^i \in \partial h(x^i)$ for~$i\ge1$.
Then Assumption~\ref{assumption:stepsize} implies~$x^k,y^k \in x^0 + \Span\{g^0,\dots,g^{k-1}\} + \Span\{s^1,\dots,s^k\}$ for~$k\ge1$.

\paragraph{Semidefinite lifting}
On top of the smooth Gram representation~\eqref{eq:gram-matrix}, we encode the prox-subgradients
given by~$s^k \in \partial h(x^k)$, the optimal subgradient~$s^\star \in \partial h(x^\star)$, and the values~$h^k = h(x^k)$.
Define the matrix~$P \in \reals^{d \times (2K+3)}$ and the vector~$F \in \reals^{2K+2}$ by
\begin{align*}
    P &= \begin{bmatrix}
        (x^0 - x^\star) & \bigl(\nabla f(y^k)\bigr)_{k=0}^{K-1} & \nabla f(x^K) & s^\star & \bigl(s^k\bigr)_{k=1}^{K}
    \end{bmatrix}, \\
    F &= \Bigl(
        \bigl(f(y^k) - (f+h)^\star\bigr)_{k=0}^{K-1}, f(x^K) - (f+h)^\star,
        h^\star,
        \bigl(h(x^k) - (f+h)^\star\bigr)_{k=1}^{K}
    \Bigr),
\end{align*}
along with the extended Gram matrix~$G = P^T P \in \symm_+^{2K+3}$.
Note that we have~$f^\star + h^\star = (f+h)^\star = 0$ from Assumption~\ref{asm:composite_h}.
When~$h \equiv 0$, the prox-subgradient columns~$s^\star, s^1, \dots, s^K$ and the~$h$-value entries~$h^\star, h(x^1) - (f+h)^\star, \dots, h(x^K) - (f+h)^\star$ disappear, recovering the conventional PEP representation~\eqref{eq:gram-matrix}.
 
\section{Proof of Proposition~\ref{prop:ls_rate}}
\label{appendix:ls_rate}
The proof has three components.
First, Lemmas~\ref{lem:spectral_moment} and~\ref{lem:accel_spectral_moment} estimate the second and fourth spectral moments of the residuals.
Next, Lemma~\ref{lem:cvar_expectation} shows that, at fixed dimension and tail level, the value-at-risk is exponentially smaller than the mean.
The final proof combines these facts to establish the CVaR equivalence in part~(i) and the second-moment bound in part~(ii).
We consider two first-order methods: gradient descent (GD) and Nesterov's fast gradient method (FGM), both with step size~$1/L$ and, for FGM, momentum coefficient~$\gamma_k = k/(k+3)$.
For either method, define the residual polynomial~$\rho_K$ by
\begin{equation*}
    x^K - x^\star = \rho_K(Q)(x^0 - x^\star).
\end{equation*}
Thus, in an eigendirection of~$Q$ with eigenvalue~$\lambda$, the initial error is multiplied by~$\rho_K(\lambda)$ after~$K$ iterations.
The residual polynomial of GD is~$\rho_K^{\text{GD}}(\lambda) = (1 - \lambda/L)^K$.
For FGM, set~$u=\lambda/L$ and~$\psi_k(u)=\rho_k^{\text{FGM}}(Lu)$.
The updates give~$\psi_0(u)=1$, $\psi_1(u)=1-u$, and
\begin{equation*}
    \psi_{k+1}(u)
    =(1-u)\left(
        \frac{2k+1}{k+2}\psi_k(u)
        -\frac{k-1}{k+2}\psi_{k-1}(u)
    \right),
    \qquad k\geq1.
\end{equation*}
\citet[Appendix~A.1.2]{paquetteHaltingTimePredictable2023} analyze this three-term recurrence through Legendre polynomials and obtain the uniform Bessel approximation in their Corollary~A.1:
\begin{equation*}
    \rho_K^{\text{FGM}} (\lambda) \sim
\tfrac{2 J_1 (K \sqrt{\lambda/L})}{K \sqrt{\lambda/L}}\, e^{-\lambda K/(2L)},
\end{equation*}
where~$J_1$ is the Bessel function of the first kind of order~$1$.

\begin{lemma}
\label{lem:spectral_moment}
    Adopt the spectral density assumption of Proposition~\ref{prop:ls_rate}, with~$\nu$ supported on~$[0, L]$.
    Then for integers~$m \ge 1$ and~$n \ge 1$, in the regime~$K \to \infty$,
    \begin{equation*}
        I_{m,n}(K) := \int_0^L \lambda^m \big(\rho_K^{\text{GD}}(\lambda)\big)^{n} p(\lambda)\, d\lambda
        \;\sim\;
        p_0\, \Gamma(m + a + 1)\, \left( \frac{nK}{L} \right)^{-(m + a + 1)}.
    \end{equation*}
\end{lemma}

\begin{proof}
    Parameterize~$\lambda$ by~$t = -\log(1 - \lambda/L)$, so~$\lambda = L(1 - e^{-t})$ and~$d\lambda = L e^{-t}\, dt$, mapping~$[0, L)$ to~$[0, \infty)$.
    This puts~$I_{m,n}(K)$ in exact Laplace form,
    \begin{equation*}
        I_{m,n}(K) = {\textstyle \int_0^\infty} e^{-nKt}\, g_m(t)\, dt,
        \quad g_m(t) := L^{m+1} (1 - e^{-t})^m e^{-t}\, p \big( L(1 - e^{-t}) \big),
    \end{equation*}
    where~$g_m$ is integrable on~$(0, \infty)$.
    From~$1 - e^{-t} \sim t$, $e^{-t} \to 1$, and~$p(\lambda) \sim p_0 \lambda^a$ as~$\lambda \to 0^+$, we obtain~$g_m(t) \sim p_0 L^{m+a+1}\, t^{m+a}$ as~$t \to 0^+$, with exponent~$m + a > -1$.
    Watson's lemma~\citep[Theorem~3.1, Chapter~3]{olverAsymptoticsSpecialFunctions1997} then gives the stated asymptotic.
\end{proof}

Unlike gradient descent, the fast gradient method has an oscillatory residual polynomial, so its second and fourth moments follow separate asymptotics from those of Lemma~\ref{lem:spectral_moment}.

\begin{lemma}
\label{lem:accel_spectral_moment}
    Adopt the spectral density assumption of Proposition~\ref{prop:ls_rate}, with~$\nu$ supported on~$[0, L]$, and let~$\rho_K^{\text{FGM}}$ be the residual polynomial of Nesterov's fast gradient method, with the Bessel asymptotic~$\rho_K^{\text{FGM}}(\lambda) \sim 2 (K\sqrt{\lambda/L})^{-1} J_1(K\sqrt{\lambda/L})\, e^{-\lambda K/(2L)}$ \citep[Corollary~A.1]{paquetteHaltingTimePredictable2023}.
    Then, as~$K \to \infty$,
    \begin{align*}
        H_2(K) &:= \int_0^L \lambda \rho_K^{\text{FGM}}(\lambda)^2 p(\lambda)\, d\lambda = \Theta \big( K^{-e(a)} \big),\\
        H_4(K) &:= \int_0^L \lambda^2 \rho_K^{\text{FGM}}(\lambda)^4 p(\lambda)\, d\lambda = \Theta \big( K^{-2h(a)} \big),
    \end{align*}
    where~$e(a)=2(a+2)$ for~$a\le-1/2$ and~$e(a)=a+7/2$ for~$a\ge-1/2$, while~$h(a)=a+3$ for~$a\le0$ and~$h(a)=a/2+3$ for~$a\ge0$.
    The rates have additional~$\log K$ multipliers at~$a=-1/2$ and~$a=0$, respectively.
\end{lemma}

\begin{proof}
    Write~$u=\lambda/L$ and
    \begin{equation*}
        b_K(u)=\frac{2e^{-Ku/2}J_1(K\sqrt{u})}{K\sqrt{u}},
        \qquad
        \delta_K(u)=\big|\psi_K(u)-b_K(u)\big|.
    \end{equation*}
    By~\citet[Corollary~A.1, equation~(108)]{paquetteHaltingTimePredictable2023}, for~$0\le u\le \log^2(K)/K$,
    \begin{equation*}
        \delta_K(u) \le C e^{-Ku/2}
        \begin{cases}
            K^{-2/3}, & 0\le u\le K^{-4/3},\\
            u^{-1/2}K^{-4/3}, & K^{-4/3}<u\le \log^2(K)/K.
        \end{cases}
    \end{equation*}
    Near the origin, the density assumption implies~$p(Lu)\le C u^a$.
    Using~$|\psi_K(u)^2-b_K(u)^2|\le C(\delta_K(u) |b_K(u)|+\delta_K(u)^2)$, $|\psi_K(u)^4-b_K(u)^4|\le C(\delta_K(u) |b_K(u)|^3+\delta_K(u)^4)$, and~$|J_1(x)|\le C\min\{x,x^{-1/2}\}$, split the weighted error integrals at~$K^{-2}$, $K^{-4/3}$, and~$\log^2(K)/K$.
    Direct integration on these intervals gives
    \begin{equation*}
        {\textstyle \int_0^1} u^{s/2}p(Lu)\big|\rho_K^{\text{FGM}}(Lu)^s-b_K(u)^s\big|\,\mathrm du
        =O\big((1+\log K)K^{-g_s(a)}\big),
        \quad \text{for } s\in \{2,4\}.
    \end{equation*}
    Here~$g_2(a)=4a/3+23/6$ for~$-1<a\le-3/4$ and~$g_2(a)=a+43/12$ for~$a\ge-3/4$, while
    \begin{equation*}
        g_4(a)=
        \begin{cases}
            2a+20/3, & -1<a\le-3/4,\\
            4a/3+37/6, & -3/4\le a\le-1/4,\\
            a+73/12, & a\ge-1/4.
        \end{cases}
    \end{equation*}
    For example, on~$[0,K^{-2}]$ we have~$b_K(u)=O(1)$ and~$\delta_K(u)=O(K^{-2/3})$, so the fourth-power error is~$O(K^{-2/3}\int_0^{K^{-2}}u^{a+2}\,\mathrm du)=O(K^{-(2a+20/3)})$.
    At the error-bound breakpoints~$a=-3/4$ and~$a=-1/4$, either adjacent expression gives the same value; the logarithmic factor is needed only there.
    The exact residual representation in~\citet[equation~(105)]{paquetteHaltingTimePredictable2023} makes the remaining tail superpolynomially small: if~$u\ge\log^2(K)/K$, then~$(1-u)^{K/2}\le e^{-\log^2(K)/2}$.
    Since~$g_2(a)-e(a)\ge1/12$ and~$g_4(a)-2h(a)\ge1/12$, both errors are smaller than their respective leading terms by a polynomial factor.

    It remains to analyze the leading Bessel terms.
    For every~$\epsilon>0$, choose~$\delta>0$ such that~$|p(Lu)/(p_0L^a u^a)-1|\le\epsilon$ for~$0<u\le\delta$.
    On~$[\delta,1]$, both the exact residual and its Bessel approximation are exponentially small in~$K$, so their contributions are negligible relative to the polynomial leading terms.
    Thus, after first letting~$K\to\infty$ and then~$\epsilon\to0$, we may replace~$p(Lu)$ by~$p_0 L^a u^a$ in the leading Bessel integrals.
    Rescaling by~$t=K^2u$ gives
    \begin{equation*}
        H_s(K)\sim 2^s p_0 L^{a+1+s/2}K^{-(2a+2+s)}
        {\textstyle \int_0^{K^2}} J_1(\sqrt t)^s e^{-st/(2K)}t^a\,\mathrm dt,
        \quad \text{for } s\in \{2,4\}.
    \end{equation*}
    With~$x=\sqrt t$, the integral factor becomes~$A_{s,K}(a):=2\int_0^K J_1(x)^s e^{-sx^2/(2K)}x^{2a+1}\,\mathrm dx$.
    Split the integral at~$x=1$, since the two relevant asymptotics of~$J_1$ hold on different ranges.
    On~$[0,1]$, we have~$J_1(x)\sim x/2$, so this block contributes a positive constant independent of~$K$.
    On~$[1,K]$, we have~$J_1(x)^s=(2/(\pi x))^{s/2}\cos^s(x-3\pi/4)+O(x^{-s/2-1})$ for~$s\in\{2,4\}$, and the exponential factor gives an effective cutoff at~$x=O(\sqrt K)$.
    Partitioning~$[1,\sqrt K]$ into fixed-length periods gives matching upper and lower bounds because~$\cos^2$ and~$\cos^4$ have positive period averages and the remaining weight varies by bounded factors on each period.
    Consequently,
    \begin{equation*}
        A_{2,K}(a) = \Theta\Big(1+{\textstyle \int_1^{\sqrt K}}x^{2a}\,\mathrm dx\Big),
        \qquad
        A_{4,K}(a) = \Theta\Big(1+{\textstyle \int_1^{\sqrt K}}x^{2a-1}\,\mathrm dx\Big).
    \end{equation*}
    The first integral changes behavior at~$a=-1/2$, the second at~$a=0$.
    The first integral factor is~$\Theta(1)$ for~$a<-1/2$, $\Theta(\log K)$ for~$a=-1/2$, and~$\Theta(K^{a+1/2})$ for~$a>-1/2$.
    The second integral factor is~$\Theta(1)$ for~$a<0$, $\Theta(\log K)$ for~$a=0$, and~$\Theta(K^a)$ for~$a>0$.
    Combining these leading terms with the smaller approximation errors proves the stated rates.
\end{proof}

Since the transfer from the average-case rate~$\Expect(\perfmetr[K])$ to~$\cvar_\alpha (\perfmetr[K])$ is identical for both GD and FGM, we isolate it in a separate lemma, which separates typical instances from the rare instances that determine the mean.
At fixed~$d$ and~$\alpha$, with probability at least~$1-\alpha$ all eigenvalues stay above a fixed positive threshold and the initial error has bounded norm.
The loss then decays exponentially.
The mean nevertheless decays only polynomially because it receives contributions from increasingly rare instances with eigenvalues near zero.
This separation makes the value-at-risk negligible relative to the mean.

\begin{lemma}
    \label{lem:cvar_expectation}
    Let~$\perfmetr[K] = \tfrac12 \sum_{i=1}^d c_i^2\, \lambda_i\, \rho_K(\lambda_i)^2$, where~$\rho_K$ is the residual polynomial of a first-order method at iteration~$K$ and the coefficients~$c_i$ satisfy~$\Expect\sum_i c_i^2=r^2<\infty$.
    Let~$\lambda_1, \dots, \lambda_d$ be the eigenvalues of a matrix~$Q \in \symm_+^{d}$, lying in~$[0, L]$ and exchangeable with common marginal~$\nu$ satisfying~$\nu(\{0\})=0$, equivalently~$\nu([0,s])\to0$ as~$s\to0^+$.
    Suppose the residual decays exponentially away from the origin: for every~$s > 0$, there is~$\kappa_s>0$ such that~$\sup_{\lambda \in [s, L]} \rho_K(\lambda)^2 = O(e^{-\kappa_s K})$.
    Also suppose that~$\Expect(\perfmetr[K]) = K^{-q+o(1)}$ for some~$q > 0$.
    Then for every tail level~$\alpha \in (0, 1]$, the relation~$\cvar_\alpha(\perfmetr[K]) \sim \alpha^{-1}\, \Expect(\perfmetr[K])$ holds.
\end{lemma}

\begin{proof}
    The case~$\alpha=1$ follows directly from~$\cvar_1(\perfmetr[K])=\Expect(\perfmetr[K])$.
    Fix~$\alpha \in (0,1)$.
    Evaluating the~$\cvar_\alpha$ representation at the optimal threshold~$\var_\alpha(\perfmetr[K])$ yields
    \begin{equation*}
        \alpha^{-1} \big( \Expect(\perfmetr[K]) - \var_\alpha(\perfmetr[K]) \big)
        \;\le\;
        \cvar_\alpha(\perfmetr[K])
        \;\le\;
        \alpha^{-1}\, \Expect(\perfmetr[K]),
    \end{equation*}
    so it suffices to show~$\var_\alpha(\perfmetr[K]) = o(\Expect(\perfmetr[K]))$.
    From~$\sum_i c_i^2 \lambda_i \le L \|x^0 - x^\star\|_2^2$ and~$\lambda_i \ge \lambda_{\min}(Q)$, the loss obeys~$\perfmetr[K] \le \tfrac12 L \|x^0 - x^\star\|_2^2 \cdot \sup_{\lambda \ge \lambda_{\min}(Q)} \rho_K(\lambda)^2$.
    Since~$\Prob(\lambda_{\min}(Q) < s) \le d\, \nu([0, s])$ and~$\nu([0, s]) \to 0$, choose~$s > 0$ so that this probability is below~$\alpha/2$.
    Markov's inequality and~$\Expect\|x^0-x^\star\|^2=r^2$ give~$\Prob(\|x^0-x^\star\|^2>2r^2/\alpha)\le\alpha/2$.
    Thus, with probability at least~$1-\alpha$, both~$\lambda_{\min}(Q)\ge s$ and~$\|x^0-x^\star\|^2\le2r^2/\alpha$ hold.
    On this event, the decay hypothesis gives~$\perfmetr[K]\le Lr^2\alpha^{-1}\sup_{\lambda\ge s}\rho_K(\lambda)^2=O(e^{-\kappa_s K})$.
    Hence the~$(1 - \alpha)$-quantile obeys~$\var_\alpha(\perfmetr[K]) = O(e^{-\kappa_s K}) = o(\Expect(\perfmetr[K]))$, which proves the claim.
\end{proof}

\begin{proof}[Proof of Proposition~\ref{prop:ls_rate}]
Writing~$x^0 - x^\star = \sum_i c_i v_i$ in the eigenbasis~$\{v_i\}_{i=1}^d$ of~$Q$, independence from~$Q$ and rotational invariance give~$\Expect(c_i^2\mid Q) = r^2/d$.
    For either method, the suboptimality~$\perfmetr[K]$ can be written as
\begin{equation*}
        \perfmetr[K](f, x^0) = \tfrac{1}{2} {\textstyle \sum_{i=1}^d} c_i^2\, \lambda_i\, \rho_K(\lambda_i)^2,
        \qquad
        \Expect(\perfmetr[K]) = \tfrac{r^2}{2} {\textstyle \int_0^L} \lambda\, \rho_K(\lambda)^2\, p(\lambda)\, d\lambda,
    \end{equation*}
    where~$\rho_K = \rho_K^{\text{GD}}$ or~$\rho_K^{\text{FGM}}$.
    The expectation identity follows directly by conditioning on~$Q$:
    \begin{equation*}
        \Expect\big(\perfmetr[K]\mid Q\big)
        =\frac{r^2}{2d}\sum_{i=1}^d \lambda_i \rho_K(\lambda_i)^2.
    \end{equation*}
    Taking expectation over~$Q$ and using exchangeability of the eigenvalues yields the fixed-density integral displayed above.

Lemma~\ref{lem:spectral_moment} with~$(m,n)=(1,2)$ gives the stated GD expectation rate, and Lemma~\ref{lem:accel_spectral_moment} gives the FGM expectation rate.
    Both residuals decay exponentially away from the origin.
    For GD, monotonicity gives~$\sup_{\lambda \in [s, L]} \rho_K^{\text{GD}}(\lambda)^2 = (1 - s/L)^{2K} = e^{-\Theta(K)}$ for every~$s > 0$.
    For FGM, the exact representation~\citep[equations~(105)--(107)]{paquetteHaltingTimePredictable2023} writes the residual at~$u=\lambda/L$ as~$2(1-u)^{(K+1)/2}I_K(u)/(K\sqrt u)$ with~$I_K$ uniformly bounded.
    Hence its supremum over~$\lambda\in[s,L]$ is~$O(e^{-\kappa_s K})$ for some~$\kappa_s>0$.
    Lemma~\ref{lem:cvar_expectation} therefore yields~$\cvar_\alpha(\perfmetr[K]) \sim \alpha^{-1}\, \Expect(\perfmetr[K])$ for both methods.

    We now prove part~(ii).
    The dual representation of~$\cvar_\alpha$ and Cauchy--Schwarz give
    \begin{align*}
        \cvar_\alpha(\perfmetr[K])
        &=
        \inf_{t \in \reals}
        \sup_{0 \le \zeta \le \alpha^{-1}}
        \left\{
            \Expect\big(\zeta\,\perfmetr[K]\big)
            + t \big(1 - \Expect\zeta\big)
        \right\}
        \\&=
        \sup_{\substack{0 \le \zeta \le \alpha^{-1}\\ \Expect\zeta=1}}
        \Expect\big(\zeta\,\perfmetr[K]\big)
        \le
        \sqrt{\alpha^{-1} \Expect\big((\perfmetr[K])^2\big)}.
    \end{align*}
    As~$x^0 - x^\star$ is isotropic Gaussian and independent of~$Q$, the coefficients are independent of the eigenvalues and satisfy~$c_i \overset{\text{i.i.d.}}{\sim} \mathcal{N}(0, r^2/d)$ conditionally on~$Q$.
    Hence~$\Expect(c_i^4) = 3 (r^2/d)^2$ and~$\Expect(c_i^2 c_j^2) = (r^2/d)^2$ for~$i \ne j$.
    From~$\Expect(w_i w_j) \le \Expect(w_i)\Expect(w_j)$ for~$i \ne j$,
    \begin{align*}
        \Expect\big( (\perfmetr[K])^2 \big)
        &=
            \frac{1}{4} \sum_i \Expect(c_i^4)\, \Expect( w_i^2)
            + \frac{1}{4} \sum_{i \ne j} \Expect(c_i^2 w_i \cdot c_j^2 w_j)
        \\&\le
            \frac{3}{4} \left( \frac{r^2}{d} \right)^2 \sum_i \Expect( w_i^2)
            + \big( \frac{1}{2} \sum_i \Expect (c_i^2 w_i) \big)^2,
    \end{align*}
    The expansion can be summarized as
    \begin{equation*}
        \Expect\big((\perfmetr[K])^2\big)
        \leq
        \frac{3r^4}{4d}
        \int_0^L \lambda^2 \rho_K(\lambda)^4 p(\lambda)\, d\lambda
        +\big(\Expect(\perfmetr[K])\big)^2.
    \end{equation*}
    The first term is the diagonal Gaussian contribution and carries the factor~$d^{-1}$; the cross-moment assumption bounds all off-diagonal contributions by the squared mean.
    For GD, the spectral integral equals~$I_{2,4}(K) = \Theta(K^{-(a+3)})$ by Lemma~\ref{lem:spectral_moment} with~$(m, n) = (2, 4)$.
    For FGM, it equals~$H_4(K) = \Theta(K^{-2 h(a)})$ by Lemma~\ref{lem:accel_spectral_moment}.
    Consequently, with~$h(a)$ defined for the corresponding method in Proposition~\ref{prop:ls_rate},
    \begin{equation*}
        \Expect\big((\perfmetr[K])^2\big)
        = O\big(d^{-1} K^{-2 h(a)}\big) + \big(\Expect(\perfmetr[K])\big)^2,
    \end{equation*}
    For FGM, the additional~$\log K$ factor at~$a=0$ is the one stated in part~(ii).
    Taking the square root gives the two-term bound in the proposition.
    For every fixed~$d$, the expectation term is lower order for all~$a>-1$; for example, the GD exponents satisfy~$a+3<2(a+2)$.
    This proves the fixed-dimensional rate, whose asymptotic threshold may depend on~$d$.
\end{proof}

\section{Proof of Theorem~\ref{thm:trajectory_implication}}
\label{appendix:trajectory_implication}

We prove the composite generalization (Theorem~\ref{thm:trajectory_implication_composite}) of Theorem~\ref{thm:trajectory_implication}, in which~$\algo$ is a fixed-step proximal gradient method applied to~$f + h$ with~$(f, h, x^0) \in (\fclass, \fclassprox)_\xclass$, \ie, $f \in \fclass$, $h \in \fclassprox$, and $\|x^0 - x^\star\| \le r$.
Theorem~\ref{thm:trajectory_implication} is recovered by setting~$h \equiv 0$, in which case~$\prox_{\gamma h}$ is the identity, the forward half-steps coincide with the first-order updates, and the gradient update coefficients need not satisfy~$\sum_{i=0}^k \eta_k^i = \gamma^{k+1}$.
The optimal subgradient~$s^\star$ vanishes, the~$h$ evaluations and prox-subgradient terms disappear, and the instance distance~\eqref{eq:composite_instance_norm} reduces to~$\|f_1 - f_2\| + \|x_1^0 - x_2^0\|$.

\begin{theorem}\label{thm:trajectory_implication_composite}
    Let~$K\ge1$ be an integer.
    Assume that the optimal subgradients are uniformly bounded, \ie, $S := \sup_{(f, h, x^0) \in (\fclass, \fclassprox)_\xclass} \|s^\star\| < \infty$.
    For~$i=1,2$, let~$\{x_i^k\}_{k=0,\dots,K}$ be generated by applying the fixed-step proximal gradient method~$\algo$ to~$f_i+h_i$, where~$f_i \in \fclass$ and~$h_i \in \fclassprox$.
    There exists~$C_K > 0$ such that
    \begin{equation*}
        \norm{Z_1 - Z_2} \le
            C_K\, \norm{(f_1, h_1, x_1^0) - (f_2, h_2, x_2^0)},
    \end{equation*}
    for any~$(f_1, h_1, x_1^0), (f_2, h_2, x_2^0) \in (\fclass, \fclassprox)_{\xclass}$, where~$\domain \subset \reals^d$ is the common convex compact domain from Lemma~\ref{lem:domain_size}.
\end{theorem}

We first show that all iterates and forward half-steps generated within the iteration budget~$K$ stay inside a common convex compact set~$\domain \subset \reals^d$.
Thanks to the normalization~$x^\star = 0$, this set is independent of the instance.

\begin{lemma}\label{lem:domain_size}
    Suppose that the integer~$K \ge 1$ is fixed, $\|x^0 - x^\star\| \le r$ for some~$r > 0$, and the optimal subgradient is uniformly bounded as~$\|s^\star\| = \|\nabla f(x^\star)\| \le S < \infty$.
    There exists~$D > 0$ such that for every instance~$(f, h, x^0) \in (\fclass, \fclassprox)_{\xclass}$, the trajectory iterates~$\{x^i\}_{i=0}^K$, the extrapolated points~$\{y^i\}_{i=0}^{K}$, and the forward half-steps~$\{x^{i+1/2}\}_{i=0}^{K-1}$ of the fixed-step proximal gradient method~$\algo$ all lie in the convex compact domain~$\domain = \{\, x \in \reals^d \mid \|x\| \le D \,\}$.
    For~$h \equiv 0$, as~$\prox_{\gamma h} = \identity$ and~$S = 0$, we recover the smooth statement used in Theorem~\ref{thm:trajectory_implication}.
\end{lemma}

\begin{proof}
    We prove by induction on~$k$ that~$\|x^i - x^\star\| \le D_k$ for all~$i \le k \le K$, where~$\bar{\eta} = \max_{k, i} |\eta_k^i|$, $B = \max_{0 \le k \le K} \sum_{j=0}^{k} |\beta_k^j| \ge 1$, and~$D_k = \big( {\textstyle \prod_{i=1}^{k}} B ( 1 + i\, \bar{\eta}\, L ) \big) r$.
    The base case~$k=0$ holds since~$\|x^0 - x^\star\| \le r = D_0$.
    As the extrapolation coefficients sum to one, the induction hypothesis gives~$\|y^i - x^\star\| \le {\textstyle \sum_{j=0}^{i}} |\beta_i^j|\, {\|x^j - x^\star\|} \le B D_k$ for~$i \le k$.
    The optimality condition gives~$x^\star = \prox_{\gamma^{k+1} h}(x^\star - \gamma^{k+1} \nabla f(x^\star))$.
    The nonexpansiveness of~$\prox_{\gamma^{k+1} h}$~\citep[Prop.~12.28]{bauschkeConvexAnalysisMonotone2017} yields
    \begin{align*}
        \|x^{k+1} - x^\star\|
        &\le
            \big\| (x^{k+1/2} - x^\star) + \gamma^{k+1} \nabla f(x^\star) \big\|
        \\&=
            \big\| (y^k - x^\star) - {\textstyle\sum_{i=0}^{k}} \eta_k^i \big( \nabla f(y^i) - \nabla f(x^\star) \big) \big\|.
    \end{align*}
    For~$h \equiv 0$, the equality in the second line holds without the coupling condition on~$\eta_k^i$ because~$\nabla f(x^\star)=0$.
    From the~$L$-smoothness of~$f$, $\|\nabla f(y^i) - \nabla f(x^\star)\| \le L \|y^i - x^\star\|$, so this bound gives~$\|x^{k+1} - x^\star\| \le B \big( 1 + (k+1) \bar{\eta}\, L \big) D_k = D_{k+1}$, hence~$\|x^i - x^\star\| \le D_K$ and~$\|y^i - x^\star\| \le B D_K$ for all~$i \le K$.
    The half-steps then satisfy~$\|x^{k+1/2} - x^\star\| \le \|y^k - x^\star\| + \sum_{i=0}^{k} |\eta_k^i|\, \|\nabla f(y^i)\| \le B D_K + \bar\eta\,(K+1)\,(S + L B D_K) =: D$, using~$\|\nabla f(y^i)\| \le S + L B D_K$.
    As~$B \ge 1$, we have~$D \ge B D_K \ge D_K$, so all three families of points lie in~$\domain$.
\end{proof}

We then bound the difference between two distinct trajectories generated from different problem instances and initial iterates.
\begin{lemma}\label{lem:trajectory_distance}
    Let~$\domain \subseteq \{x\in\reals^d\mid\|x\|\le D\}$ be compact and convex, with~$0\in\domain$.
    Suppose that~$f_1,f_2$ are continuously differentiable on~$\domain$, have~$L$-Lipschitz gradients, and satisfy~$\|\nabla f_i(0)\|\le S$ for~$i=1,2$.
    Let~$\varepsilon_f = \sup_{x \in \domain} |f_1(x) - f_2(x)|$ and~$\varepsilon_g = \sup_{x \in \domain} \|\nabla f_1(x) - \nabla f_2(x)\|$.
    For any~$x_1, x_2 \in \domain$ with~$\varepsilon_x = \| x_1 - x_2 \|$, the following bounds hold:
    \begin{equation*}
        \big| f_1(x_1) - f_2(x_2) \big|
        \le
            (S + L D)\, \varepsilon_x
            + \varepsilon_f,
        \qquad
        \big\| \nabla f_1(x_1) - \nabla f_2(x_2) \big\|
        \le
            L\, \varepsilon_x
                + \varepsilon_g.
    \end{equation*}
\end{lemma}

\begin{proof}
    For~$i=1,2$ and~$x \in \domain$, we have~$\|\nabla f_i(x)\| \le \|\nabla f_i(0)\| + \|\nabla f_i(x) - \nabla f_i(0)\| \le S + L D$.
    Using the mean value theorem,~$|f_1(x_1) - f_2(x_2)| \le |f_1(x_1) - f_1(x_2)| + |f_1(x_2) - f_2(x_2)| \le (S + LD)\,\varepsilon_x + \varepsilon_f$.
    Furthermore, $\|\nabla f_1 (x_1) - \nabla f_2 (x_2)\| \le \|\nabla f_1 (x_1) - \nabla f_1 (x_2)\| + \|\nabla f_1 (x_2) - \nabla f_2 (x_2)\| \le L\,\varepsilon_x + \varepsilon_g$.
\end{proof}

Before proving the theorem, we define the distance metric over~$(\fclass, \fclassprox)_{\xclass}$ using Moreau envelopes for the nonsmooth part.

\begin{definition}
    \label{def:composite_instance_norm}
    For~$\gamma>0$, the Moreau envelope~$M_{\gamma h}(v) := \min_x\{ h(x) + (2\gamma)^{-1}\|x - v\|^2 \}$ of~$h \colon \reals^d \to \reals \cup \{\infty\}$ is finite-valued and continuously differentiable if~$h$ is proper, lower semicontinuous, and convex.
    Write the~$C^1$ norm~$\|\phi\|_{C^1(\domain)} := \sup_{x \in \domain} |\phi(x)| + \sup_{x \in \domain} \|\nabla \phi(x)\|$ over the compact set~$\domain \subset \reals^d$.
    We identify functions whose values and gradients agree on~$\domain$.
    Define the distance $\| (f_1, h_1, x_1^0) - (f_2, h_2,x_2^0) \|$ between~$(f_1, h_1, x_1^0)$ and~$(f_2, h_2, x_2^0)$ in~$(\fclass, \fclassprox)_{\xclass}$ by
    \begin{equation}\label{eq:composite_instance_norm}
\norm[C^1(\domain)]{f_1 - f_2}
            + \max_{1 \le k \le K}\! \norm[C^1(\domain)]{M_{\gamma^k h_1} \!- M_{\gamma^k h_2}}
            + \norm{x_1^0 - x_2^0}.
    \end{equation}
\end{definition}

We now prove Theorem~\ref{thm:trajectory_implication_composite},
which reduces to the smooth convex case when~$h \equiv 0$.

\begin{proof}[Proof of Theorem~\ref{thm:trajectory_implication_composite}]
Define
\begin{equation*}
        \varepsilon_x=\|x_1^0-x_2^0\|,\qquad
        \varepsilon_f=\sup_{v\in\domain}|f_1(v)-f_2(v)|,\qquad
        \varepsilon_g=\sup_{v\in\domain}\|\nabla f_1(v)-\nabla f_2(v)\|,
    \end{equation*}
    and
    \begin{equation*}
        \varepsilon_M = \max_{1 \le k \le K} \|M_{\gamma^k h_1} - M_{\gamma^k h_2}\|_{C^1(\domain)},
        \quad
        \varepsilon_{\mathrm{prox}} = 
            \max_{1 \le k \le K}\,
\sup_{v \in \domain}\,
            \bigl\|\prox_{\gamma^k h_1}(v) - \prox_{\gamma^k h_2}(v)\bigr\|.
    \end{equation*}
Write~$\bar\varepsilon = \varepsilon_x + \varepsilon_f + \varepsilon_g + \varepsilon_{\mathrm{prox}} + \varepsilon_M$ and~$\delta = \|(f_1, h_1, x_1^0) - (f_2, h_2, x_2^0)\|$.
    Since $\nabla M_{\gamma^k h}(v) = (\gamma^k)^{-1}\big(v - \prox_{\gamma^k h}(v)\big)$, we have~$\varepsilon_{\mathrm{prox}} \le \bar\gamma\, \varepsilon_M$ with~$\bar\gamma = \max_k \gamma^k$, from~$\prox_{\gamma^k h_1}(v) - \prox_{\gamma^k h_2}(v) = -\gamma^k \big(\nabla M_{\gamma^k h_1}(v) - \nabla M_{\gamma^k h_2}(v)\big)$.
    As each~$\varepsilon_x$, $\varepsilon_f$, $\varepsilon_g$, and~$\varepsilon_M$ is at most~$\delta$ by definition~\eqref{eq:composite_instance_norm}, we have~$\bar\varepsilon \le c_0\, \delta$ for some~$c_0 > 0$.
    Also write~$\underline\gamma = \min_k \gamma^k > 0$, $\bar{\eta} = \max_{k, i} |\eta_k^i|$, and~$B = \max_k \sum_j |\beta_k^j| \ge 1$ as in Lemma~\ref{lem:domain_size}, noting that the two trajectories share these algorithm coefficients.
    First, we have
    \begin{equation}\label{eq:prox-assembly-decomp}
        \norm{Z_1 - Z_2}^2
        \le \norm[2]{P_1 + P_2}^2\, \norm[F]{P_1 - P_2}^2 + \norm[2]{F_1 - F_2}^2.
    \end{equation}
    The columns of~$P_i$ are~$x_i^0 - x^\star$, the gradients~$\nabla f_i(y_i^k)$ for~$k<K$, the terminal gradient~$\nabla f_i(x_i^K)$, $s_i^\star = -\nabla f_i(x^\star)$, and the prox-subgradients~$s_i^k$.
    For any~$v\in\domain$, nonexpansiveness and~$x^\star=\prox_{\gamma^k h_i}(x^\star+\gamma^k s_i^\star)$ give
    \begin{equation*}
        \|\prox_{\gamma^k h_i}(v)-x^\star\|
        \le \|v-x^\star-\gamma^k s_i^\star\|
        \le D+\gamma^kS.
    \end{equation*}
    Hence~$\|\nabla M_{\gamma^k h_i}(v)\|\le2D/\gamma^k+S\le2D/\underline\gamma+S=:S_{\mathrm{prox}}$.
    In particular,~$s_i^k=\nabla M_{\gamma^k h_i}(x_i^{k-1/2})$ has norm at most~$S_{\mathrm{prox}}$.
    Each column of~$P_1 + P_2$ thus has norm at most~$2 \max\{ D, S + LD, S_{\mathrm{prox}} \}$, so~$\|P_1 + P_2\|_{2} \le \|P_1 + P_2\|_F \le 2 \sqrt{2K+3}\, \max\{ D, S + LD, S_{\mathrm{prox}} \} =: C$.

    Set~$\Delta_k = \|x_1^k - x_2^k\|$, $\Delta_k^y = \|y_1^k - y_2^k\|$, $\Delta_{k+1/2} = \|x_1^{k+1/2} - x_2^{k+1/2}\|$, and~$\bar\Delta_k = \max_{0 \le j \le k} \Delta_j$.
    The nonexpansiveness of the proximal map and Lemma~\ref{lem:trajectory_distance} give
    \begin{align*}
        \Delta_k^y &\le {\textstyle\sum_{j=0}^k} |\beta_k^j|\,\Delta_j \le B\bar\Delta_k, \\
        \Delta_{k+1/2}
        &\le \Delta_k^y + \bar\eta {\textstyle\sum_{j=0}^k} (L\Delta_j^y + \varepsilon_g)
        \le \chi\bar\Delta_k + (K+1)\bar\eta\varepsilon_g, \\
        \Delta_{k+1} &\le \Delta_{k+1/2} + \varepsilon_{\mathrm{prox}},
    \end{align*}
    where~$\chi = B(1 + (K+1)\bar\eta L) \ge 1$.
    Thus~$\bar\Delta_{k+1} \le \chi\bar\Delta_k + c\bar\varepsilon$ for~$c = 1 + (K+1)\bar\eta$, so induction from~$\bar\Delta_0 = \varepsilon_x \le \bar\varepsilon$ bounds~$\Delta_k$, $\Delta_k^y$, and~$\Delta_{k+1/2}$ by~$\hat D_K\bar\varepsilon$ for some finite~$\hat D_K$.
    Another application of Lemma~\ref{lem:trajectory_distance} at the points~$y_i^k$ for~$k<K$ and~$x_i^K$ bounds the corresponding gradient and function-value gaps by finite multiples of~$\bar\varepsilon$; also,~$\|s_1^\star-s_2^\star\| \le \varepsilon_g$.

    It remains to bound the~$h$-value entries and prox-subgradient columns.
First, for~$k = 0, \dots, K-1$, the Moreau decomposition~$h_i(x_i^{k+1}) = M_{\gamma^{k+1} h_i}(x_i^{k+1/2}) - \tfrac{\gamma^{k+1}}{2} \|s_i^{k+1}\|^2$ gives
\begin{equation*}
        \big| M_{\gamma^{k+1} h_1}(x_1^{k+1/2}) - M_{\gamma^{k+1} h_2}(x_2^{k+1/2}) \big|
        \le \varepsilon_M + S_{\mathrm{prox}}\, \|x_1^{k+1/2} - x_2^{k+1/2}\|,
    \end{equation*}
    where splitting~$\varepsilon_M$ uses the function-value part of the Moreau~$C^1$ norm, not the prox-distance alone, as the envelopes need not agree at any reference point.
    Second, for~$k=1,\dots,K$,
    \begin{equation*}
        \norm{s_1^k - s_2^k}
        =
            \norm{(x_1^{k-1/2} - x_1^k)/\gamma^{k} - (x_2^{k-1/2} - x_2^k)/\gamma^{k}}
        \le
            \underline{\gamma}^{-1} \left( \Delta_{k-1/2} + \Delta_k \right),
    \end{equation*}
    Together with~$|\|s_1^k\|^2 - \|s_2^k\|^2| \le 2 S_{\mathrm{prox}}\, \|s_1^k - s_2^k\|$, this gives~$|h_1(x_1^k) - h_2(x_2^k)| \le E_K\, \bar\varepsilon$ for a finite~$E_K$.
    Moreover,~$|h_1^\star - h_2^\star| = |f_1^\star - f_2^\star| \le \varepsilon_f$, where the equality follows from~$f_i^\star+h_i^\star=0$.

    Finally, each column difference of~$P_1 - P_2$ and each entry difference of~$F_1 - F_2$ is at most a constant multiple of~$\bar\varepsilon$.
    Substituting into~\eqref{eq:prox-assembly-decomp} with~$\|P_1 + P_2\|_{2} \le C$ and~$\bar\varepsilon \le c_0\, \delta$ yields~$\|Z_1 - Z_2\| \le C_K\, \delta$ for a finite~$C_K > 0$.
\end{proof}

\begin{remark}[Subgradient methods]\label{rem:subgradient_methods}
    The \emph{subgradient method} for nonsmooth convex minimization replaces~$g^k = \nabla f(x^k)$ with~$g^k \in \partial f(x^k)$, typically over the class of convex Lipschitz functions.
    The statement of Theorem~\ref{thm:trajectory_implication_composite} need not hold here, since the subdifferential~$\partial f$ is set-valued and different selections~$x \mapsto g \in \partial f(x)$ yield distinct trajectories.
    For instance, in dimension~$d=1$, the function~$f(x) = M|x|$ satisfies~$\partial f(0) = [-M, M]$.
    Even~$f_1 = f_2 = f$ (hence~$M_{\gamma f_1} = M_{\gamma f_2}$) admits selections~$g_i \in \partial f_i(0)$ with~$|g_1 - g_2| = 2M > 0$.
    Therefore, the trajectory gap is uncontrolled by any instance distance.
\end{remark}

\section{Proof of Proposition~\ref{prop:worst-case-expectations}}
\label{appendix:worst-case-expectations}

\begin{proof}[Proof of Proposition~\ref{prop:worst-case-expectations}]
    For the first inequality, we first show that the ambient instance space~$\fclassX$ is compact under~\eqref{eq:c1_sup_norm}.
    By Lemma~\ref{lem:domain_size}, the iterates of every instance~$(f, x^0) \in \fclassX$ lie in the fixed convex compact domain~$\domain = \{\, x \in \reals^d \mid \|x\| \le D \,\}$, and the metric~\eqref{eq:c1_sup_norm} sees only the restrictions of~$f$ and~$\nabla f$ to~$\domain$.
    The normalization~$x^\star = 0$, $f^\star = 0$ gives~$\|\nabla f(x)\| \le LD$ and~$|f(x)| \le LD^2$, while~$\nabla f$ is~$L$-Lipschitz.
    Hence the families~$\{f|_\domain \mid f \in \fclass\}$ and~$\{\nabla f|_\domain \mid f \in \fclass\}$ are uniformly bounded and uniformly equicontinuous, so the Arzel\`a--Ascoli theorem makes them relatively compact in the supremum norm.
    They are closed under joint uniform convergence because convexity, $L$-smoothness, and the normalization are preserved.
    Combining this with the compact initial-iterate ball~$\xclass(f) = \{\, x \in \reals^d \mid \|x - x^\star\| \le r \,\}$ proves that~$\fclassX$ is compact.

    The support of~$\prob$ is therefore compact.
    The empirical distribution~$\hprob$ converges weakly to~$\prob$ almost surely by the strong law for empirical measures on separable metric spaces~\citep[Theorem~11.4.1]{dudleyRealAnalysisProbability2002}; compactness upgrades this to~$\lim_{N \to \infty} W_{\fclassX}(\prob, \hprob) = 0$ almost surely~\citep[Theorem~6.9]{villaniOptimalTransportOld2009}.
    For every~$\beta \in (0,1)$, define the deterministic sequence~$\{ \varepsilon_N(\beta) \}$ for $N\ge 1$ by
    \begin{equation*}
        \varepsilon_N (\beta)
        = C_K \, \inf 
        \bigl\{
            t > 0 \;\bigm|\;
            \prob^N \big( W_{\fclassX}(\prob, \hprob) \le t \big) \ge 1-\beta
        \bigr\}.
    \end{equation*}
    This sequence is well-defined because~$\fclassX$ is compact and therefore bounded.
    Note that almost sure convergence gives~$\prob^N \big( W_{\fclassX}(\prob, \hprob) \le t \big) \to 1$ for any~$t > 0$.
    Since~$1 - \beta < 1$, this probability is at least~$1 -\beta$ for all~$N$ large enough,
    so $\varepsilon_N(\beta) \le C_K t$ for all such~$N$ and~$\varepsilon_N(\beta) \to 0$ follows.

    Explicit rates for such a radius on infinite-dimensional function spaces are available~\citep{leiConvergenceConcentrationEmpirical2020}, but are not needed here.

    For the second inequality, it suffices to show that the lifting~$(f, x^0) \mapsto Z$ maps every distribution in~$\mathcal{P}=\mathcal{B}_{\fclassX}(\hprob,\varepsilon/C_K)$ to a distribution in~$\ambiset[\varepsilon]$ with the same expected loss.
    Considering the performance metric~$\perfmetr[K]$ in Section~\ref{sec:prob_analysis_fom}, this lifting preserves the loss, \ie, $\ell(f, x^0) = \ell(Z)$, and is~$C_K$-Lipschitz by Theorem~\ref{thm:trajectory_implication}.
    Therefore, for every pair of distributions~$\probQ$ and~$\hprob$ on~$\fclassX$, their pushforwards satisfy~$W(T_\#\probQ,T_\#\hprob)\le C_K W_{\fclassX}(\probQ,\hprob)$, where~$T(f,x^0)=Z$.
    This proves the claim.
\end{proof}

\end{document}